\documentclass[11pt,reqno]{amsart}

\usepackage{amsfonts}
\usepackage{eurosym}
\usepackage{amssymb}
\usepackage{amsthm}
\usepackage{amsmath}
\usepackage{amsaddr}
\usepackage{bm}
\usepackage{cite}
\usepackage{mathrsfs}
\usepackage{xcolor}
\usepackage[OT1]{fontenc}
\usepackage[left=1.8cm, right=1.8cm, top=3cm, bottom=3cm]{geometry}
\usepackage{hyperref}
\usepackage{graphicx}
\usepackage{mathtools}
\usepackage{cleveref}
\usepackage{comment}
\hypersetup{colorlinks=true, linkcolor=blue, citecolor=red, urlcolor=blue}

\RequirePackage{times}

\usepackage{xpatch}
\makeatletter
\AtBeginDocument{\xpatchcmd{\@thm}{\thm@headpunct{.}}
{\thm@headpunct{}}{}{}}
\makeatother

\allowdisplaybreaks
\newtheorem{theorem}{Theorem}[section]

\newtheorem{definition}[theorem]{Definition}
\newtheorem{lemma}[theorem]{Lemma}

\newtheorem{remark}[theorem]{Remark}

\def \au {\rm}
\def \ti {\it}
\def \jou {\rm}
\def \bk {\it}
\def \no#1#2#3 {{\bf #1} (#3), #2.}
\def \eds#1#2#3 {#1, #2, #3.}

\def\ve{\varepsilon}

\newcommand \norm[1] {\left\Vert #1\right\Vert}
\def\n{\textbf{\textit{n}}}
\def\R{\mathbb{R}}

\def\G{\mathcal{G}}
\def\F2o{\overline{F_2}}

\def\d{{\rm d}}
\def \l {\langle}
\def \r {\rangle}

\def\ddt{\frac{\d}{\d t}}

\def\ve{\varepsilon}
\def\vphi{\varphi}

\def\div{\mathrm{div}\,}
\def\pt{\partial_t}

\def \au {\rm}
\def \ti {\it}
\def \jou {\rm}
\def \bk {\it}
\def \no#1#2#3 {{\bf #1} (#3), #2.}
\def \eds#1#2#3 {#1, #2, #3.}
\def \nome#1#2 {{\bf #1}, (#2).}

\newcommand{\sref}[2]{\hyperref[#2]{#1 \ref*{#2}}}

\numberwithin{equation}{section}

\makeatletter
\def\@settitle{\begin{center}%
		\baselineskip14\p@\relax
		\huge
		\@title
	\end{center}%
}
\makeatother

\makeatletter
\@namedef{subjclassname@2020}{%
  2020 Mathematics Subject Classification}
\makeatother

\begin{document}

\title[Unconstrained approximation of the Cahn--Hilliard model]{\emph{On the dynamics of an unconstrained phase-field approximation  of the Cahn--Hilliard equation}}

\author[M. Conti, P. Galimberti, S. Gatti, A. Giorgini, G. Ricchi]{Monica Conti$^1$, Pietro Galimberti$^2$, Stefania Gatti$^3$, Andrea Giorgini$^1$, Greta Ricchi$^1$}
\address{$^1$Politecnico di Milano\\
	Dipartimento di Matematica\\
	Via E. Bonardi 9, I-20133 Milano, Italy \\
	\href{mailto:monica.conti@polimi.it}{monica.conti@polimi.it},
	\href{mailto:andrea.giorgini@polimi.it}{andrea.giorgini@polimi.it},
	\href{mailto:greta.ricchi@polimi.it}{greta.ricchi@polimi.it}}

\address{$^2$Università di Ferrara\\
	Dipartimento di Matematica e Informatica\\
	Via Machiavelli 30, I-44121, Ferrara, Italy\\
	\href{mailto:pietro.galimberti@edu.unife.it}{pietro.galimberti@edu.unife.it}
	}

\address{$^3$Università di Modena e Reggio Emilia\\
	Dipartimento di Scienze Fisiche, Informatiche e Matematiche\\
	Via Campi 213/B, I-41125 Modena, Italy
	\\
	\href{mailto:stefania.gatti@unimore.it}{stefania.gatti@unimore.it}
	}

\keywords{Cahn--Hilliard equation, logarithmic potential, mobility, well-posedness, strict separation property.}

\subjclass[2020]{35K35, 35K65, 35B25}

\begin{abstract}
We study a generalized Cahn--Hilliard equation, based on an unconstrained theory proposed by Duda, Sarmiento and Fried, with non-constant mobility and nonlinear potential of logarithmic type in two and three space dimensions. In the case of non-degenerate mobility, we prove well-posedness of weak solutions, propagation of regularity and strict separation for the phase-field variable. Concerning the degenerate mobility case, we show existence of a class of weak solutions. Compared with the concept of weak solution for the classical Cahn--Hilliard equation introduced by Elliott and Garcke, in our case the chemical potential is identified pointwise.
Finally, we prove stability with respect to the perturbation parameter, showing that weak solutions of the generalized system converge to weak solutions of the classical Cahn--Hilliard equation on any finite time interval, in both mobility regimes.
\end{abstract}

\maketitle

\setcounter{tocdepth}{1}

\tableofcontents

\section{Introduction}
The classical Cahn--Hilliard equation reads as
\begin{alignat}{2}
	\label{CH1}
	\pt \vphi &= \div(b(\vphi) \nabla \mu) \quad && \text{in} \ \Omega \times (0, \infty),\\
	\label{CH2}
	\mu &= -\Delta\vphi + \Psi'(\vphi) && \text{in} \ \Omega \times (0, \infty),
\end{alignat}
equipped with the following classical boundary and initial conditions
\begin{equation}  \label{nCH-mu0}
	\partial_\n \vphi= b(\varphi) \partial_\n \mu=0 \quad \text{on } \, \partial \Omega \times (0,\infty),\quad
	\vphi|_{t=0}= \vphi_0 \quad \text{in } \, \Omega.
\end{equation}
This system was originally proposed by Cahn and Hilliard in \cite{CH} and \cite{CH71} as a model to describe phase separation in binary alloys. Here, the set $\Omega \subset \R^d$, $d=2,3$, is the physical domain occupied by the two components. The phase-field variable $\vphi:\Omega \times (0, \infty) \to [-1,1]$ represents the difference between the concentrations of the two materials with $-1$ and $1$ being the pure phases. The chemical potential $\mu$ is the variational derivative of the Ginzburg--Landau free energy
\begin{equation}
	\label{free energy}
	E(\vphi)= \int_{\Omega} \frac12 |\nabla \vphi |^2 + \Psi(\vphi)\,  \d x,
\end{equation}
where $\Psi$ is the double-well Flory--Huggins potential
\begin{equation}
	\label{FH}
	\Psi(s) = F(s)- \frac{\theta_0}{2}s^2 = \frac{\theta}{2}((1+s)\ln(1+s)+(1-s)\ln(1-s))- \frac{\theta_0}{2}s^2, \quad s \in [-1,1],
\end{equation}
with positive parameters $\theta$ and $\theta_0$ satisfying the relation $\theta_0-\theta>0$. The function $b : [-1,1] \to \R$, introduced in the original derivation in \cite{CH}, is the Onsager mobility and measures the strength of the diffusion. Significant choices for $b$ are
\begin{align}
	&b(s)\equiv m \in \R_+ && \text{constant case} \label{m-const}\\
	&b \in C([-1,1]):\ 0 < b_m \le  b(s) \le b_M  \quad &&\text{non degenerate case} \label{m-non deg}\\
	&b(s)= m (1-s^2), \ m \in \R_+ &&\text{degenerate case} \label{m-deg}
\end{align}

Over the years, numerous variants of the Cahn--Hilliard equation have been proposed to account for additional physical mechanisms. Coupling with fluid equations yields Navier--Stokes/Cahn--Hilliard and Hele--Shaw/Cahn--Hilliard systems for incompressible two-phase flows \cite{ABELS2009,AbelsDepnerGarcke2013,GGW2018}, while coupling with reaction--diffusion processes leads to models of nutrient transport and tumour growth \cite{ColliGilardiHilhorst2015,GarckeLam2016,GarckeLam2017}. Other versions include the Cahn--Hilliard--Oono equation \cite{Oono1987,Miranville2011,GGM}, hyperbolic models accounting for finite relaxation of the diffusion flux \cite{Galenko, GGPM, GSSZ, Segatti2007}, and second-order relaxations of the fourth-order problem \cite{PP2022,CFGP2025}; see \cite{MB2019} for a broader overview. Despite their differences, most of these models retain the constitutive structure of the classical theory and modify only its dynamics.

A different line of research questions the constitutive foundations of the Cahn--Hilliard equation itself. This perspective originates in the continuum thermodynamics of order parameters developed by Fried and Gurtin \cite{FG1993, FG1994}. In \cite{Gurtin}, Gurtin derived both the Ginzburg--Landau and the Cahn--Hilliard equations by supplementing the balance of mass with a separate balance law for microforces, namely generalized forces conjugate to the evolution of the order parameter, together with constitutive relations compatible with the second law of thermodynamics. In this setting the Cahn--Hilliard equation is no longer postulated phenomenologically, and natural generalizations, such as viscous Cahn--Hilliard-type equations \cite{Novick-Cohen}, arise from different constitutive choices.

Building on this framework, Duda, Sarmiento and Fried \cite{fried} recently proposed a theory of coupled diffusion and phase transition in which the concentration $w$ and the phase-field variable $\vphi$ are independent quantities. The classical Cahn--Hilliard model corresponds to the constrained theory in which the two variables are identified, whereas relaxing this constraint leads to a new class of models. In the unconstrained case, the system proposed in \cite{fried} reads as
\begin{alignat}{2}
	\label{w problem}
	\ve \partial_t w  &= \div(b(\vphi) \nabla w) -\div(b(\vphi) \nabla \vphi)\quad &\text{in} \ \Omega \times (0,\infty),\\
	\label{w problem 2}
	w&= -\ve \Delta \vphi + \ve \Psi'(\vphi) + \vphi & \text{in} \ \Omega \times (0,\infty),
\end{alignat}
completed with boundary and initial conditions
\begin{equation}
	\label{w problem bc i condition}
	\partial_\n \vphi= b(\varphi) \partial_\n w=0 \quad \text{on } \, \partial \Omega \times (0,\infty), \quad w|_{t=0}= w_0 \quad \text{in} \ \Omega.
\end{equation}
Here $\ve>0$ measures the strength of the coupling between $w$ and $\vphi$. Indeed, the associated free energy to \eqref{w problem}--\eqref{w problem 2} is given by
$$
E_\ve(\vphi, w)=\int_\Omega \frac12|\nabla\vphi|^2+\Psi(\vphi)\,\d x+\frac{1}{2\ve}\int_\Omega|w-\vphi|^2\,\d x.
$$
The last term in $E_\ve$ penalizes deviations from the constraint $w=\vphi$, so that the classical Cahn--Hilliard equation is formally recovered as $\ve\to0$.

System \eqref{w problem}--\eqref{w problem 2} subject to \eqref{w problem bc i condition} can be rewritten as perturbation of system \eqref{CH1}--\eqref{CH2} subject to \eqref{nCH-mu0}. 
Define
\begin{equation}\label{change of variable}
\mu=\frac{w-\vphi}{\varepsilon}.
\end{equation}
If the initial data are related by $w_0=-\varepsilon\Delta\vphi_0+\varepsilon\Psi'(\vphi_0)+\vphi_0$, then system \eqref{w problem}--\eqref{w problem 2} is equivalent to
\begin{alignat}{2}
	\label{CH1-p}
	\partial_t \vphi &= \div( b(\vphi) \nabla \mu) - \ve \partial_t \mu \quad &&\text{in }  \Omega \times (0,\infty),\\
	\label{CH2-p}
	\mu &= -\Delta \vphi + \Psi^{\prime}(\vphi)  \quad &&\text{in } \Omega \times (0,\infty),
\end{alignat}
completed with boundary and initial conditions
\begin{equation}  \label{nCH-mu}
	\partial_\n \vphi= b(\varphi) \partial_\n \mu=0 \quad \text{on } \, \partial \Omega \times (0,\infty),\quad
	\vphi|_{t=0}= \vphi_0 \quad \text{in } \, \Omega.
\end{equation}

In contrast with the classical Cahn--Hilliard equation \eqref{CH1}--\eqref{CH2}, the spatial average of the phase field
\begin{equation*}
	\overline{\vphi}(t):= \frac{1}{|\Omega|} \int_{\Omega} \vphi(x,t) \,\d x
\end{equation*}
is in general not conserved along solutions to \eqref{CH1-p}--\eqref{CH2-p} subject to \eqref{nCH-mu}. What is conserved is the spatial average of the concentration $w$, namely
\begin{equation}
	\label{conservation}
	\overline{w}(t)=\frac{1}{|\Omega|}\int_{\Omega} w(x,t)\, \d x = \frac{1}{|\Omega|}\int_{\Omega} w_0(x) \,\d x, \quad \forall\, t \ge 0.
\end{equation}
Recalling that $w=\vphi+\ve\mu$, this gives
$$
\overline{\vphi}(t)=\overline{\vphi_0}+\ve\big(\overline{\mu_0}-\overline{\mu}(t)\big),
$$
so the mean value of $\vphi$ varies only through terms of order $\ve$. This is consistent with the interpretation of the classical (constrained) theory as the limit, as $\ve\to0^+$, of the unconstrained one. Indeed, the authors of \cite{fried} suggested that ``the Cahn--Hilliard equation should be recovered in the limit as $\ve$ tends to $0$'', that is, solutions to \eqref{CH1-p}--\eqref{CH2-p} converge to solutions to \eqref{CH1}--\eqref{CH2}. 

System \eqref{CH1-p}--\eqref{CH2-p} was first studied in \cite{Miranville 2021}, where it is regarded as an Allen--Cahn-type equation for an unconstrained order parameter, with constant mobility and a polynomial double-well potential. There, the author proved well-posedness and established the convergence to the Cahn--Hilliard equation as $\ve\to0$, together with an error bound. These results were extended in \cite{CGM2024} to singular potentials such as \eqref{FH}, still with constant mobility. In addition, the authors proved that $\vphi$ satisfies the strict separation property, also in three dimensions: for every $\tau>0$ there exists $\delta=\delta(\tau, \varepsilon)>0$ such that
$$
|\vphi(x, t)| \le 1- \delta \quad \forall\, (x, t) \in \overline{\Omega} \times [\tau, \infty).
$$
This result was obtained by exploiting the regularity of strong solutions. The case of degenerate mobility was considered in \cite{SLYM}, where existence of weak solutions was proved for a polynomial potential, which does not enforce the constraint $|\vphi|\le1$. We also mention that the relaxation term $\ve\pt\mu$ appearing in \eqref{CH1-p}, combined with an additional viscous term $\beta\pt\vphi$ in the chemical potential \eqref{CH2-p}, arises in several Cahn--Hilliard models of tumor growth coupled with nutrient and reaction equations \cite{ColliGilardiRoccaSprekels2015, ColliGilardiRoccaSprekels2017, Signori2020, HilhorstKampmannNguyenVanDerZee2015, ColliGilardiHilhorst2015}. In that context, well-posedness, singular limits as the relaxation parameters vanish (including quantitative error estimates) and optimal control problems have been investigated. 

To the best of our knowledge, the case of a singular potential combined with a non-constant mobility has not been addressed so far. In this paper we fill this gap by studying \eqref{w problem}--\eqref{w problem 2} subject to \eqref{w problem bc i condition}, or equivalently \eqref{CH1-p}--\eqref{CH2-p} subject to \eqref{nCH-mu}, with a singular potential such as \eqref{FH} and both non-degenerate and degenerate mobilities, in two and three space dimensions. Our contributions are twofold.

\medskip

\textbf{Non-degenerate mobility.} We extend the results of \cite{CGM2024} from constant to non-degenerate mobility. We first prove global-in-time existence of weak solutions satisfying an energy equality, together with their uniqueness when $b\in C^1([-1,1])$, and then establish propagation of regularity, with particular emphasis on the strict separation from the pure phases. This is not a mere adaptation of \cite{CGM2024}. When $b$ depends on $\vphi$, the inverse operator $\mathcal{G}_{\vphi}$ of the elliptic operator $-\div(b(\vphi)\nabla\,\cdot\,)$ depends on the solution itself. As a consequence, the stability estimate for the difference of two solutions $\vphi_1$ and $\vphi_2$ contains additional terms involving $b'(\vphi_1)\pt\vphi_1$ and $\mathcal{G}_{\vphi_1}-\mathcal{G}_{\vphi_2}$ (cf.~\cite{CGGG2025} for the classical Cahn--Hilliard equation in two dimensions). We point out that similar terms reappear in the error estimate for the singular limit.

Moreover, in contrast with \cite{CGM2024}, our proof of the separation property relies only on the regularity of weak solutions. The key idea is to exploit the parabolic structure of \eqref{w problem}: by means of a De Giorgi iteration scheme (see \cite{DiBenedettoDegenerate}), we obtain an $L^\infty(\Omega\times[\tau,\infty))$ bound for $w=\vphi+\ve\mu$, for every $\tau>0$. The elliptic problem \eqref{w problem 2} then yields the separation of $\vphi$ from the pure phases for all positive times, also in three dimensions.

Next, we prove that, as $\ve\to0^+$, weak solutions to the perturbed system \eqref{CH1-p}--\eqref{CH2-p} converge, up to a subsequence, to a weak solution to the classical Cahn--Hilliard equation \eqref{CH1}--\eqref{CH2}, in both two and three dimensions. The proof relies only on $\ve$-uniform estimates derived from the energy identity, rather than on the estimates for strong solutions used in \cite{CGM2024}. In two dimensions, we also establish the rate of convergence
$$
\|\vphi^\ve(t)-\vphi^0(t)\|_{(H^1(\Omega))'}^2+\int_0^t\|\vphi^\ve(s)-\vphi^0(s)\|_{H^1(\Omega)}^2\,\d s\le C\ve,\qquad t\in[0,T].
$$
To this aim, we perform higher-order estimates that are uniform with respect to the parameter $\ve$. The main obstacle is the nonlinear term $P$ in \eqref{P} that can not be controlled as in \cite{CGGG2025} due to the perturbation term $\ve \partial_t \mu$ appearing in \eqref{CH1-p}.
We overcome this difficulty by means of a generalized H\"older inequality from \cite{GGW2022} (see Lemma \ref{L2 product lemma} below), combined with a logarithmic Gr\"onwall argument.

\medskip

\textbf{Degenerate mobility.} For the classical Cahn--Hilliard equation with degenerate mobility, the main difficulty lies in the identification of the chemical potential $\mu$, which is not controlled when $\vphi$ reaches the pure phases. Hence, only a weaker notion of solution is available (see Def. \ref{def-CH-weak-deg} below), formulated in terms of the flux $\mathbf{J}$ rather than of $\mu$ (see \cite{EG1996}, and also \cite{AbelsDepnerGarcke2013, CFGP2025, DD2016, LMS2012}). For the perturbed system \eqref{CH1-p}--\eqref{CH2-p} with the physically relevant choice $b=\theta/F''$, which gives $b(s)=1-s^2$ for the logarithmic potential \eqref{FH}, a better structure can be exploited.
For every fixed $\ve>0$, the energy identity contains the term $\frac{\ve}{2}\|\mu\|^2_{L^2(\Omega)}$, which provides a bound on $\mu$ in $L^\infty(0,T;L^2(\Omega))$. This allows us to identify both the chemical potential $\mu$ and the derivative of the potential $\Psi'(\vphi)$ almost everywhere. As a consequence, we show that $|\vphi(x,t)|<1$ for almost every $x\in\Omega$ and for all $t>0$, whereas only $|\vphi|\le 1$ almost everywhere in $\Omega\times(0,\infty)$ is known for the classical Cahn-Hilliard equation \cite{EG1996}.

Finally, we study the asymptotic as $\ve\to0^+$. We show that, up to a subsequence, weak solutions to \eqref{CH1-p}--\eqref{CH2-p} converge to a weak solution, in the sense of \cite{EG1996}, to the classical Cahn--Hilliard equation with degenerate mobility. The key step is a bound on $\vphi^\ve$ in $L^2(0,T;H^2(\Omega))$ that is uniform with respect to $\ve$. In the literature, this is usually obtained through a differential inequality for the entropy density $\Phi$, which compensates the degeneracy of the mobility (see, e.g.,  \cite[Lemma 2]{EG1996} and also \cite{AbelsDepnerGarcke2013, CFGP2025, DD2016, LMS2012}). However, in the perturbed system, the same argument yields the additional term
$\varepsilon\langle\partial_t\mu^\ve,\,\Phi'(\varphi^\ve)\rangle$,
which has no sign and it is not a time derivative of a controlled quantity.
Hence, the entropy method can not be applied in this case.
To overcome this issue, we propose a different approach in Subsection \ref{sub-conv-deg}, which is of independent interest: we derive an $L^2(0,T;H^2(\Omega))$ estimate directly from the dissipation term in the energy identity, without using the entropy function. In fact, expanding $b(\vphi^\ve) |\nabla\mu^\ve|^2$ with $\mu^\ve=-\Delta\varphi^\ve+\Psi'(\varphi^\ve)$ gives (see \eqref{conto smart})
$$
b(\vphi^\ve)|\nabla\mu^\ve|^2=b(\vphi^\ve)|\nabla\Delta\varphi^\ve|^2+b(\vphi^\ve)(\Psi''(\vphi^\ve))^2|\nabla\varphi^\ve|^2
-2\,b(\vphi^\ve)\Psi''(\vphi^\ve)\,\nabla\Delta\varphi^\ve\cdot\nabla\varphi^\ve.
$$
Our argument hinges on the structure of the cross term: since $bF''=\theta$, after integration in time-space, we obtain the desired norm of $\Delta\varphi^\ve$ in $L^2(\Omega \times (0,T))$ as in \eqref{stima-deg} below. 

\section{Main results}

We first state the standing assumptions on $\Omega$, $\Psi$, and $\ve$ used throughout the paper.	

\subsection*{Assumptions} $\,$ 

\medskip\noindent
$[\text{H1}]\,$
The set $\Omega \subset \R^d$, $d=2,3$ is a bounded domain with boundary $\partial \Omega$ of class $C^3$.

\medskip\noindent
$[\text{H2}]\,$ The nonlinearity $\Psi$ is a quadratic perturbation of
a singular (strictly) convex function in $[-1, 1]$:
$$\Psi(s)=F(s)-\frac{\theta_0}{2}s^2,$$
where the convex part $F$ belongs to $C([-1,1])\cap C^{2}(-1,1)$ and fulfills
\begin{equation*}
	\label{ipoH}
	\lim_{s\rightarrow -1}F^{\prime }(s)=-\infty ,\quad \lim_{s\rightarrow
		1}F^{\prime }(s)=+\infty,
	\quad F^{\prime \prime }(s)\geq
	\theta, \quad\forall\, s\in (-1,1),
\end{equation*}
for some $\theta>0$. Here, we study the physical case of double-well potentials, namely, we assume
$
\theta_0>\theta.
$
We also extend $F(s)=+\infty$ for any $s\notin [-1,1]$.
Without loss of generality, we assume that
$F(0)=F'(0)=0.$ In particular, this entails that $F(s)\geq 0$ for all $s\in [-1,1]$ and $F'(s)s\geq 0$. 
The assumptions are motivated by the logarithmic potential defined in \eqref{FH}.

\medskip
\noindent
$[\text{H3}]\,$ We assume that $\ve\in(0, \ve_0)$, where $\ve_0>0$ satisfies
\begin{equation}
	\label{epps}
	\ve_0 < \frac{1}{\theta_0}.
\end{equation}
Condition \eqref{epps} guarantees that the map $s \mapsto \ve\Psi'(s)+s$  is monotone (see, e.g., \eqref{phi H1 w L2} and \eqref{norm phi e norm w} below).
\medskip

\subsection*{Non-degenerate mobility} 
The first part of the paper is devoted to the analysis of system \eqref{CH1-p}--\eqref{CH2-p}, subject to \eqref{nCH-mu}, in the case of non-degenerate mobility $b$ satisfying the assumptions in \eqref{m-non deg}. We begin by establishing a well-posedness result. The notation used throughout the paper is collected in Section \ref{Math Setting}.

\begin{theorem}
\label{ws thm phi and mu}
Let $\ve\in (0,\ve_0)$ be given. Assume that the initial datum $\vphi_0 \in H^2(\Omega)$ is such that $\norm{\vphi_0}_{L^{\infty}(\Omega)}\le1$, $|\overline{\vphi_0}|<1$, $\partial_\n \vphi_0=0$ on $\partial \Omega$, and
$$\mu_0:= -\Delta\vphi_0+ \Psi'(\vphi_0) \in L^2(\Omega).$$
	\begin{itemize}
	\item[(A)] \textbf{Existence and uniqueness of weak solutions.}  There exists a global weak solution $\vphi$ to \eqref{CH1-p}--\eqref{nCH-mu} in the following sense:
\begin{align*}
&\vphi \in L^{\infty}(0, \infty; H^2(\Omega)), \quad \pt \vphi \in L^2(0, \infty; H^1(\Omega)),
\nonumber\\
&\mu \in L^{\infty}(0, \infty; L^2(\Omega)) \cap L^2_{\rm uloc}([0, \infty); H^1(\Omega)), \quad \pt \mu \in L^2(0, \infty; (H^1(\Omega))'), 
\nonumber \\
&F'(\vphi) \in L^{\infty}(0, \infty; L^2(\Omega)),
\end{align*}
with
$$
	| \vphi(x,t)| <1 \quad \text{a.e. } x \in \Omega, \ \forall \, t \ge 0,
$$
		such that
		\begin{equation}
		\label{weak formulation perturbed Cahn--Hilliard equation}
		(\pt \vphi, v) + \ve \l\pt \mu, v \r + (b(\vphi)\nabla \mu, \nabla v)=0 \quad \quad \forall \, v \in H^1(\Omega),
		\end{equation}
		for almost every $t \in (0,\infty)$, and
\begin{equation}
\label{mu perturbed Cahn--Hilliard equation}
\mu =  -\Delta \vphi + \Psi'(\vphi) \quad \text{a.e. in } \Omega \times (0, \infty).
\end{equation}
Furthermore, $\partial_{\n} \vphi =0$ a.e. on $\partial \Omega \times (0, \infty)$, $\vphi(0)=\vphi_0$ a.e. in $\Omega$ and the energy equality
\begin{align*}
&\frac12 \norm{\nabla\vphi (t)}^2_{L^2(\Omega)}
+ \frac{\ve}{2}\norm{\mu(t)}_{L^2(\Omega)}^2 
+ \int_\Omega \Psi(\vphi(t) ) \, \d x 
+ \int_0^t \int_\Omega b(\vphi)|\nabla \mu |^2 \, \d x \, \d s
\\
& \quad =\frac12 \norm{\nabla\vphi_0}^2_{L^2(\Omega)}
+ \frac{\ve}{2}\norm{\mu_0}_{L^2(\Omega)}^2 
+ \int_\Omega \Psi(\vphi_0) \, \d x
\end{align*}
holds for all $t \in (0,\infty)$.
Suppose, in addition, that $b \in C^1([-1,1])$. Then, the weak solution is unique.
\medskip
		
\item[(B)] \textbf{$L^\infty$-regularity and separation from the pure phases.} For any $\tau >0$, there exists $C=C(\ve,\tau)$ such that
\begin{equation*}
\norm{\vphi+\ve \mu}_{L^{\infty}( \Omega \times [\tau, \infty))} \le C.
\end{equation*}
Moreover, there exists $\delta=\delta(\ve, \tau) >0$ such that
\begin{equation*}
|\vphi(x, t)| \le 1 - \delta \quad \forall (x, t)\in \overline{\Omega}\times [\tau, \infty).
\end{equation*}

\item[(C)] \textbf{Strong solutions.} Assume, in addition, that $\mu_0 \in H^1(\Omega)$ and $b \in C^1([-1,1])$. Then, the unique solution $\vphi$ is a strong solution in the following sense:
\begin{align*}
&\vphi \in L^{\infty}(0, \infty; W^{2,p}(\Omega)), \quad \pt \vphi \in L^\infty(0, \infty; H^1(\Omega)), \nonumber\\
&\mu \in L^{\infty}(0,\infty; H^1(\Omega)) \cap L^2_{\rm uloc}([0, \infty); H^2(\Omega)), \quad \pt \mu \in L^2_{\rm uloc}([0, \infty); L^2(\Omega)),\nonumber\\
&F'(\vphi)  \in L^{\infty}(0, \infty; L^{p}(\Omega)),
\end{align*}
where $p=6$ if $d=3$ and any $2 \le p <\infty$ if $d=2$.
\end{itemize}
\end{theorem}

Next, we focus on the convergence of solutions to \eqref{CH1-p}--\eqref{CH2-p} towards solutions to \eqref{CH1}--\eqref{CH2} as the perturbation parameter $\varepsilon\to 0$. We first recall the notion of weak solution to the limit problem $(\varepsilon= 0)$ with non-degenerate mobility satisfying \eqref{m-non deg} (see, e.g., \cite{BB1999, GP2025}).
\begin{definition}
\label{def-CH-weak}
Let $\vphi_0 \in H^1(\Omega)\cap L^\infty(\Omega)$ be such that $\norm{\vphi_0}_{L^{\infty}(\Omega)} \le 1$ and $|\overline{\vphi_0}| < 1$. A weak solution $\vphi: \Omega \times [0,\infty) \to \mathbb{R}$ to \eqref{CH1}--\eqref{CH2}, subject to \eqref{nCH-mu0}, is such that
\begin{align*}
&\vphi \in BC([0,\infty); H^1(\Omega)) \cap L^4_{\rm uloc}([0, \infty); H^2(\Omega)) \cap L^2_{\rm uloc}([0, \infty); W^{2,p}(\Omega)), 
\\
&\vphi \in L^{\infty}(\Omega \times (0, \infty)): \ |\vphi(x,t)| < 1 \quad \text{a.e. in }\ \Omega \times (0, \infty),
\\
&\partial_t \vphi \in L^2(0, \infty; (H^{1}(\Omega))'), \quad F'(\varphi)\in L_{\rm uloc}^2([0,\infty); L^6(\Omega)), 
\\
&\mu \in L_{\rm uloc}^2([0,\infty); H^1(\Omega)),
\end{align*} 
where $p=6$ if $d=3$ and any $2 \le p <\infty$ if $d=2$, and
\begin{align*}
&\l\partial_t \vphi, v \r + ( b(\vphi) \nabla \mu, \nabla v) = 0,  \quad \forall \, v \in H^1(\Omega), \ \text{for a.e.} \, t \in (0,\infty),
\\
&\mu= - \Delta \vphi + \Psi'(\vphi) \quad \text{a.e. in } \, \Omega \times (0,\infty).
\end{align*} 
Moreover, $\partial_\n \varphi=0$ almost everywhere on $\partial \Omega \times (0,\infty)$, and $\varphi|_{t=0}=\varphi_0$.
\end{definition}
We can now state the following result.
\begin{theorem}
\label{th-singular-limit}
Let the assumptions of Theorem \ref{ws thm phi and mu} be in place, and let $T>0$ be fixed. 
Let $\lbrace \vphi^\ve\rbrace_{\ve \in(0,\ve_0)}$ be a family of weak solutions to the perturbed problem \eqref{CH1-p}--\eqref{CH2-p}, subject to \eqref{nCH-mu} with the same initial condition $\vphi_0$.
\begin{itemize}
\item[(A)]
There exists a sub-sequence, still labelled by $\vphi^\ve$, such that
\begin{alignat*}{2}
\vphi^\ve &\to \vphi^\ast \quad &&\text{strongly in } C([0,T]; L^2(\Omega))\cap L^4(0, T; H^{2- \alpha}(\Omega)),\ \forall \, \alpha \in(0,2],
\\
\varepsilon \mu^\varepsilon &\to 0 \quad &&\text{strongly in } C([0, T]; L^2(\Omega)),
\\
b(\vphi^\ve)\nabla \mu^\ve &\rightharpoonup b(\vphi^\ast) \nabla \mu^\ast \quad &&\text{weakly in}\ L^2(0, T; L^2(\Omega;\mathbb{R}^d)),
\nonumber
\end{alignat*}
as $\varepsilon\to 0$, where $\vphi^\ast$ is a weak solution to the Cahn--Hilliard equation \eqref{CH1}--\eqref{CH2} in the sense of Definition \ref{def-CH-weak} with the same initial datum. 
\item[(B)]	Let $d=2$. 
Suppose that $b\in C^1([-1,1])$. Then, under the additional hypothesis $\mu_0 \in H^1(\Omega)$, there exist $\ve_1\in (0,\ve_0)$ and $C>0$ such that the following error bound holds
$$
\| \vphi^\ve (t)-\vphi^\ast(t)\|_{(H^{1}(\Omega))'}^2+\int_0^T\|\vphi^\ve(\tau)-\vphi^\ast(\tau)\|_{H^1(\Omega)}^2 	\, \d \tau \leq C \ve,
$$
for any $t\in [0,T]$ and $\ve\in(0,\ve_1)$, where $C>0$ depends on $T$ and on the norm of the initial datum.		
\end{itemize}
\end{theorem}

\subsection*{Degenerate mobility} 
The second part of the paper focuses on the perturbed model \eqref{CH1-p}--\eqref{CH2-p}, subject to \eqref{nCH-mu}, when the mobility $b$ is vanishing at the pure phases $\pm 1$. Precisely, we assume that the mobility $b \in C^1([-1,1])$ is given by
\begin{equation}
\label{mobility and potential compatibility condition}
b(s)= \frac{\theta}{F''(s)}, \quad \forall \, s \in (-1,1), \qquad b(\pm1)=0.
\end{equation}
Notice that \eqref{mobility and potential compatibility condition} and $F''\ge\theta$ entail $0<b\le 1$ in $(-1,1)$, while $b(\pm1)=0$ since $F''(s)\to+\infty$ as $s\to\pm1$.

A main issue when dealing with degenerate mobility is that the notion of weak solutions is very mild, the identification of the chemical potential being a major obstacle, see \cite{EG1996}. On the other hand, for the model under study such an identification is possible. More precisely, our result reads as follows.
\begin{theorem}
\label{th-degenerate}
Let the initial datum $\vphi_0 \in H^2(\Omega)$ be with $\norm{\vphi_0}_{L^{\infty}(\Omega)}\le1, |\overline{\vphi_0}|<1$, and $\partial_\n \vphi_0=0$ on $\partial \Omega$, such that
$$\mu_0:= -\Delta\vphi_0+ \Psi'(\vphi_0) \in L^2(\Omega).$$
Then, there exists a global weak solution $\vphi$ to \eqref{CH1-p}--\eqref{CH2-p} such that, for all $T>0$,
\begin{alignat*}{2}
&\vphi \in L^{\infty}(0, T; H^2(\Omega)),\quad &&\pt \vphi \in L^2(0, T; H^1(\Omega)),\nonumber\\
&\mu \in L^{\infty}(0, T; L^2(\Omega)), &&\pt \mu \in L^2(0, T; (H^1(\Omega))'), \nonumber \\
&F'(\vphi) \in L^{\infty}(0, T; L^2(\Omega)),
\end{alignat*}
with
$$
 | \vphi(x,t)| <1 \quad \text{a.e. } x \in \Omega, \ \forall \, t \ge 0,
$$
and the following relations hold:
\begin{equation}
\label{weak formulation perturbed Cahn--Hilliard equation-deg}
(\pt \vphi, \eta) + \ve \langle \pt\mu,\eta \rangle - (\mu, \div(b(\vphi)\nabla \eta))=0,
\end{equation}
for any $\eta \in H^2(\Omega)$  with $\partial_{\n}\eta=0$ on $\partial\Omega$, for almost every $t \in (0,\infty)$, and 
$$ 
\mu =  -\Delta \vphi + \Psi'(\vphi)\quad \text{a.e. in } \Omega \times (0, \infty),
$$
Furthermore, $\mu$ satisfies
\begin{equation}
\label{link mu J}
\left(\mu, \div(b(\vphi) \bm{\xi}) \right)
=
-\left( \Delta \vphi, \div(b(\vphi) \bm{\xi})\right)
- \left( (b\Psi'')(\vphi)\nabla\vphi, \bm{\xi} \right)\!,
\end{equation}
for any $\bm{\xi} \in H^1(\Omega; \mathbb{R}^d)$ such that $\bm{\xi}\cdot \n=0$ on $\partial\Omega$, for almost every $t \in (0,\infty)$.
In addition, there exists $\mathbf J\in L^2(\Omega\times(0,\infty); \mathbb{R}^d)$ 
such that 
\begin{equation}
\label{J-eps-0}
(\mathbf J,\bm\xi)=-(\mu,\div(b(\vphi)\bm\xi))
\end{equation}
for all $\bm{\xi} \in H^1(\Omega; \mathbb{R}^d)$ such that $\bm{\xi}\cdot \n=0$ on $\partial\Omega$ and for almost every $t\in(0,\infty)$. Finally, $\partial_{\n} \vphi =0$ a.e. on $\partial \Omega \times (0, \infty)$ and $\vphi(0)=\vphi_0$ in $\Omega$.
\end{theorem}

Lastly, we study the asymptotics as $\varepsilon\to 0$ of the weak solutions to \eqref{CH1-p}--\eqref{CH2-p} with degenerate mobility.
Let us recall the notion of weak solution  to the limit problem $(\varepsilon=0)$ as in \cite{EG1996}.
\begin{definition}
\label{def-CH-weak-deg}
Let $\vphi_0 \in H^1(\Omega)\cap L^\infty(\Omega)$ be such that $\norm{\vphi_0}_{L^{\infty}(\Omega)} \le 1$ and $|\overline{\vphi_0}| < 1$. 
A weak solution to \eqref{CH1}--\eqref{CH2} with degenerate mobility, subject to \eqref{nCH-mu0}, is a pair 
$(\vphi,\mathbf J)$ such that
\begin{align*}
&\vphi \in L^\infty(0,\infty; H^1(\Omega)) \cap L_{\rm uloc}^2([0,\infty); H^2(\Omega)), \quad \partial_t \vphi \in L^2(0, \infty; (H^{1}(\Omega))'),
\\
&\vphi \in L^{\infty}(\Omega \times (0, \infty)): \ |\vphi(x,t)| \leq  1 \quad \text{a.e. in }\ \Omega \times (0, \infty),
\\
&\mathbf J\in L^2(\Omega\times(0,\infty); \mathbb{R}^d),
\end{align*} 
which satisfies, for any $T>0$,
\begin{equation*}
\int_0^T\l\pt\vphi,v\r\,\d s+\int_0^T\int_\Omega\mathbf J\cdot\nabla v\,\d x\,\d s=0\quad\forall\,v\in L^2(0,T;H^1(\Omega)),
\end{equation*}
and 
\begin{equation*}
\int_0^T\int_\Omega\mathbf J\cdot\bm\xi\,\d x\,\d s=\int_0^T\int_\Omega\Delta\vphi\,\div(b(\vphi)\bm\xi)+(b\Psi'')(\vphi)\nabla\vphi\cdot\bm\xi\,\d x\,\d s
\end{equation*}
for all $\bm\xi\in L^2(0,T;H^1(\Omega; \mathbb{R}^d))\cap L^\infty(\Omega\times(0,T); \mathbb{R}^d)$ with $\bm\xi\cdot\n=0$ on $\partial\Omega\times(0,T)$, as well as $\partial_\n \varphi=0$ on $\partial \Omega$ almost everywhere on $(0,\infty)$, and $\varphi|_{t=0}=\varphi_0$.
\end{definition}

\begin{remark}\label{rem-EG}
We observe that the solution given by Theorem \ref{th-degenerate} is also a weak solution in the sense of Elliott and Garcke (adapted to the perturbed problem), i.e. a couple $\left(\vphi, \mathbf{J}\right)$ such that
	\begin{equation*}
	\int_0^T \l \pt (\vphi + \ve \mu), v \r  \, \d s
	+ \int_0^T \int_\Omega \mathbf{J} \cdot \nabla v \, \d x \, \d s =0,
\end{equation*}
for any $v\in L^2(0,T; H^1(\Omega))$, with $\mathbf{J}$ satisfying 
\eqref{J-eps-0}.
\end{remark}

We are now in a position to state our final result.
\begin{theorem}
\label{th-degenerate-conv}
Let the assumptions of Theorem \ref{th-degenerate} be in place. Then, the family of weak solutions $(\vphi^\ve, \mathbf{J}^\ve)_{\ve \in(0,\ve_0)}$ constructed in Theorem \ref{th-degenerate} converges, up to a sub-sequence, in the following sense:
\begin{alignat*}{2}
\vphi^\ve &\to \vphi^\ast \quad &&\text{strongly in } C([0,T]; L^2(\Omega))\cap L^2(0, T; H^{2- \alpha}(\Omega)),\ \forall \, \alpha \in(0,2], \\
\varepsilon \mu^\varepsilon &\to 0 \quad &&\text{strongly in } L^2(0, T; L^2(\Omega)),\\
		\mathbf{J}^\ve &\rightharpoonup \mathbf{J}^\ast && \text{weakly in } L^2(0, T; L^2(\Omega; \mathbb{R}^d)), \nonumber
	\end{alignat*}
as $\varepsilon \to 0$, where $\left(\vphi^\ast, \mathbf{J}^\ast \right)$ is a weak solution to the Cahn--Hilliard equation \eqref{CH1}--\eqref{CH2} with degenerate mobility, in the sense of Definition \ref{def-CH-weak-deg}, with the same initial datum.
\end{theorem}
\section{Mathematical setting}
\label{Math Setting}
\subsection{Notation and function spaces}
Let $X$ be a Banach space. We denote by $\norm{\cdot}_X$ its norm, by $X'$ its dual space and by $\l \cdot, \cdot \r_{X',X}$ the associated duality product.
Given $ 1 \le p \le \infty$ and an interval $I \subset [0, \infty)$, we denote with $L^p( I; X)$ the space of all Bochner $p-$integrable functions defined from $I$ into $X$. We also define the space $W^{1,p}(I; X)$ as the set of functions $f \in  L^p( I; X)$ with distributional derivative $ \partial_t f \in L^p(I; X)$. In particular we denote $H^1(I; X) = W^{1,2}(I; X)$. Moreover, $L^p_{\rm uloc}([0,\infty); X)$ denotes the space of functions $f\in L^p(0,T; X)$ for any $T>0$ such that
$$
\sup_{t \ge 0} \int_t^{t+1}\norm{f(s)}_X^p \d s <\infty.
$$
For any $1\leq p \leq \infty$, $L^p(\Omega)$ denotes the usual Lebesgue space. For $p=2$, we denote the inner product in $L^2(\Omega)$ with $(\cdot, \cdot)$.
For any positive integer $k, \ W^{k,p}(\Omega)$ is the Sobolev space of functions in $L^p(\Omega)$ such that the distributional derivative of order up to $k$ is an element of $L^p(\Omega)$. In particular, we use the notation $H^k(\Omega)$ for the Hilbert space $W^{k,2}(\Omega)$ with norm $\norm{\cdot}_ {H^k(\Omega)}$.
We will denote the duality product $\l \cdot , \cdot \r_{H^1(\Omega)', H^1(\Omega)}$ by $\l \cdot,\cdot \r$. 
For vector-valued functions, the corresponding Sobolev spaces are denoted by $W^{k,p}(\Omega;\mathbb{R}^d)$, where $k\in\mathbb{N}$ and $1\leq p\leq\infty$. For simplicity, we use the notation $\|\cdot\|_{W^{k,p}(\Omega)}$ for their norms, as in the scalar case.

We recall the Hilbert interpolation inequality
\begin{equation}
\label{Hilbert interpolation inequality}
\norm{f}_{L^2(\Omega)}
\le
\norm{f}_{H^1(\Omega)}^{\frac{1}{2}}
\norm{f}_{(H^1(\Omega))'}^{\frac{1}{2}}, \quad \forall \, f \in H^1(\Omega).
\end{equation}
Extending the definition of total mass as
\begin{equation*}
\overline{f} := \frac{1}{|\Omega|} \l f,1 \r,
\end{equation*}
we report the Poincaré inequality
\begin{equation}
\label{poincare}
\norm{f}_{L^2(\Omega)}
\le
C \left( \norm{\nabla f}_{L^2(\Omega)}
+ |\overline{f}| \right), \quad \, \forall f \in H^1(\Omega),
\end{equation}
where $C$ depends only on $\Omega$. For any $m \in \R$, we introduce the function space
$$
H^1_{(m)}(\Omega)=\{f \in H^1(\Omega) \ : \ \overline{f}=m\}.
$$
We recall that if $m=0$, $H^1_{(0)}(\Omega)$ is an Hilbert space with inner product $(\nabla f, \nabla g)$.
It is also convenient to introduce the Hilbert spaces of zero-mean function
\begin{equation*}
L_{(0)}^2(\Omega)=\{f \in L^2(\Omega) \ : \ \overline{f}=0\},
\quad
H^{-1}_{(0)}(\Omega)=\{f \in (H^1(\Omega))' \ : \ \overline{f}=0 \},
\end{equation*}
with norm defined as $\norm{f}_{L^2_{(0)}(\Omega)}=\norm{f}_{L^2(\Omega)}$ and $\norm{f}_{H^{-1}_{(0)}(\Omega)} = \norm{f}_{(H^1(\Omega))'}$, respectively.

\subsection{Analytic tools}

We recall the following key inequalities from \cite{GGW2022} that we will use in our analysis.
\begin{lemma}
\label{L2 product lemma}
Let $\Omega$ be a bounded domain in $\R^2$ with boundary $\partial\Omega$ of class $C^2$. Assume that $f \in H^1(\Omega)$, $g \in L^p(\Omega)$ for some $p >2$. Then, we have
	\begin{equation*}
		\norm{f g}_{L^2(\Omega)} \le C\left(\frac{p}{p-2}\right)^\frac12 \norm{f}_{H^1(\Omega)} \norm{g}_{L^2(\Omega)} \ln^\frac12 \left(e |\Omega|^{\frac{p-2}{2p}}\frac{\norm{g}_{L^p(\Omega)}}{\norm{g}_{L^2(\Omega)}}\right),
	\end{equation*}
	for some positive constant $C>0$.
\end{lemma}
\begin{lemma}
	\label{logaritmic Gronwall}
	Let $f$ be a positive absolutely continuous function on $[0, T]$ and $g, h$ be two non-negative summable functions on $[0, T]$ that satisfy the differential inequality
	\begin{equation*}
		\ddt f(t) \le g(t) f(t) \ln(e+f(t))+ h(t),
	\end{equation*}
	for almost every $t \in  [0, T]$. Then it holds
	\begin{equation*}
		f(t) \le (e+ f(0))^{e^{\int_0^t g(\tau) \d \tau}} e^{\int_0^t e^{\int_\tau^t g(s) \d s} h(\tau) \d \tau}, \quad \forall \, t \in [0, T].
	\end{equation*}
\end{lemma}
\subsection{Elliptic problems}\label{ell-prob}
We consider the Laplace equation with homogeneous Neumann boundary conditions
\begin{equation}
\begin{cases}
	\label{NP}
-\Delta w = u \ & \text{in} \ \Omega, \\
\partial_\n w = 0 \ & \text{on} \ \partial\Omega,
\end{cases}
\end{equation}
and the associated linear operator $A:H^1(\Omega) \to (H^1(\Omega))'$ defined by
\begin{equation*}
\l Aw, v \r = \int_{\Omega} \nabla w \cdot \nabla v \, \d x, \quad \forall \, w, v \in H^1(\Omega).
\end{equation*}
Exploiting the Poincaré inequality, the restriction of $A: H^1_{(0)}(\Omega) \to H^{-1}_{(0)}(\Omega)$ is an isomorphism. We denote its inverse map by $\G: H^{-1}_{(0)}(\Omega) \to H^1_{(0)}(\Omega)$, namely  for every $u \in H^{-1}_{(0)}(\Omega), \  \G u \in H^1_{(0)}(\Omega)$ is the unique weak solution of the Neumann problem \eqref{NP} such that
\begin{equation}\label{Gdef}
(\nabla \G u, \nabla v) = \l u, v \r, \quad \forall \, v \in H^1_{(0)}(\Omega).
\end{equation}
Moreover, $\norm{\nabla \G v}_{L^2(\Omega)}$ is a norm in $H^{-1}_{(0)}(\Omega)$, which is equivalent to $\norm{v}_{H^{-1}_{(0)}(\Omega)}$ defined above.

Next, let $b \in C(\mathbb{R})$ such that $0<b_m \leq b(s) \leq b_M$ for any $s \in \mathbb{R}$. Given a measurable function $q$, we consider the following elliptic problem
\begin{equation}\label{epgq}
\begin{cases}
-\div(b(q) \nabla w) = u \quad & \text{in}\ \Omega, \\
b(q) \nabla w \cdot \n = 0 \quad & \text{on} \ \partial\Omega.
\end{cases}
\end{equation}
Similarly to the definition of $\G$ in \eqref{Gdef}, we define the weak solution operator $\G_q: H^{-1}_{(0)}(\Omega) \to H^1_{(0)}(\Omega)$ as follows: for every $u \in H^{-1}_{(0)}(\Omega)$, $\G_q u \in H^1_{(0)}(\Omega)$ is the unique solution to the variational problem
\begin{equation}
\label{def G_q}
( b(q) \nabla \G_q u ,\nabla v ) = \l u, v \r, \quad \forall \, v \in H_{(0)}^1(\Omega).
\end{equation}
We recall from \cite{BB1999} that $\norm{\sqrt{b(q)}\nabla\G_q\,\cdot\,}_{L^2(\Omega)}$ defines a norm in $H^{-1}_{(0)}(\Omega)$ which is equivalent to the usual one, namely there exists $C_b>1$ such that
\begin{equation}
\label{equivalence of norm}
\frac{1}{C_b} \norm{v}_{H^{-1}_{(0)}(\Omega)}^2 \le \norm{\sqrt{b(q) }\nabla \G_{q} v}_{L^2(\Omega)}^2 \le C_b \norm{v}_{H^{-1}_{(0)}(\Omega)}^2,
\quad \forall \, v \in H^{-1}_{(0)}(\Omega).
\end{equation}
Finally, assuming in addition that $b\in C^1(\R)$ with $b'$ bounded, we report the following bounds from \cite{BB1999} and \cite{CGGG2025}: there exists $C>0$ such that
\begin{alignat}{2}
\label{H2 norm of G_q-2d}
\norm{\G_q v}_{H^2(\Omega)} &\le  C  \norm{v}_{L^2(\Omega)} 
	+C\norm{\nabla q}_{L^2(\Omega)}\norm{q}_{H^2(\Omega)}\norm{\nabla \G_q v}_{L^2(\Omega)}, \quad &&d=2,  \\
	\label{H2 norm of G_q-3d}
	\norm{\G_q v}_{H^2(\Omega)} &\le  C  \norm{v}_{L^2(\Omega)} +C \norm{q}_{H^2(\Omega)}^2\norm{\nabla \G_q v}_{L^2(\Omega)}, \quad &&d=3,
\end{alignat}
for any $q \in H^2(\Omega)$ and $v \in L^2_{(0)}(\Omega)$.

\section{Non-degenerate mobility: well-posedness}
\label{weak solution}

In this section, we prove Theorem \ref{ws thm phi and mu}. Firstly, we observe that the two formulations \eqref{w problem}--\eqref{w problem bc i condition} and \eqref{CH1-p}--\eqref{nCH-mu} of the problem are equivalent, as shown in \cite[Remark 4.3]{CGM2024}, through the change of variable \eqref{change of variable}. In light of this equivalence  Theorem \ref{ws thm phi and mu} can be reformulated into the following Theorem \ref{ws thm}.

\begin{theorem}\label{ws thm}
Let $\ve\in (0,\ve_0)$ be given. Assume the mobility $b$ satisfies \eqref{m-non deg}. Then, the following results hold:
\begin{itemize}
\item[(A)] \textbf{Existence of weak solutions.}
Let $w_0= \vphi_0 + \ve \mu_0 \in L^2(\Omega)$ with $\overline{w_0}=m \in \R$. Then, there exists a global weak solution $(w, \vphi)$ to \eqref{w problem}--\eqref{w problem 2}, subject to \eqref{w problem bc i condition}, such that
\begin{alignat}{2}
\label{ws regularity}
&w \in L^{\infty}(0, \infty; L^2(\Omega)) \cap L^2_{\rm uloc}([0, \infty); H^1_{(m)}(\Omega)), \quad \pt w \in L^2(0, \infty; H^{-1}_{(0)}(\Omega)), \nonumber
\\
&\vphi \in L^{\infty}(0, \infty; H^2(\Omega)), \quad \pt \vphi \in L^2(0, \infty; H^1(\Omega)),
\\
&F'(\vphi) \in L^{\infty}(0, \infty; L^2(\Omega)),
\nonumber
\end{alignat}
with
$$
| \vphi(x,t)| <1 \quad \text{a.e. } x \in \Omega, \ \forall \, t \ge 0.
$$
The couple $(w,\vphi)$ satisfies
\begin{equation}
\label{w-weak}
\l \ve \pt w, v\r + (b(\vphi) \nabla w, \nabla v) = (b(\vphi) \nabla \vphi, \nabla v) \quad \forall \, v  \in H^1(\Omega),
\end{equation}
for almost every $t \in (0,\infty)$, and
 \begin{equation}
 	\label{w = phi + ve mu}
 	 w = -\ve \Delta \vphi + \ve \Psi'(\vphi) + \vphi \quad \text{a.e. in } \Omega \times (0, \infty).
 \end{equation}
Furthermore, $\partial_{\n} \vphi =0$ almost everywhere on $\partial \Omega \times (0, \infty)$ and $w(0)=w_0$ almost everywhere in $\Omega$.
\smallskip

\item[(B)] \textbf{Uniqueness of weak solutions.} Suppose $b \in C^1([-1,1])$. Let $w_1, w_2$ be two weak solutions originating from initial conditions $w_1^0, w_2^0$ as in part (A). Then, there exists a constant $C>0$, depending on the $L^2$-norm of the initial data, $\ve$ and $b$, such that
\begin{equation*}
	\norm{w_1(t)-w_2(t)}_{H^{-1}_{(0)}(\Omega)}+\norm{\vphi_1(t)-\vphi_2(t)}_{H^1(\Omega)} \le C e^{Ct} \norm{w_1^0-w_2^0}_{H^{-1}_{(0)}(\Omega)} \quad \forall \, t \geq 0.
\end{equation*}
In particular, the weak solution is unique.
\smallskip

\item[(C)]\textbf{Further regularity and separation from the pure phases.} For any $\tau >0$, there exists $C=C(\ve,\tau)$ such that
\begin{equation*}
	\norm{w}_{L^{\infty}(\Omega \times [\tau, \infty))} \le C(\ve, \tau).
\end{equation*}
Moreover, there exists $\delta=\delta(\ve, \tau) >0$ such that
\begin{equation}\label{separation property}
|\vphi(x, t)| \le 1 - \delta \quad  \forall \, (x, t)\in \overline{\Omega}\times [\tau, \infty).
\end{equation}
\end{itemize}
\end{theorem}
\begin{proof}
\label{sect-existence}
We construct weak solutions to \eqref{w problem}--\eqref{w problem bc i condition} on any $[0,T]$, $T>0$, by means of a double
approximation. We first replace the singular potential $F$ with a
regularized potential $F_\lambda$, for $\lambda>0$. For each fixed
$\lambda$, we then apply a Galerkin scheme to the resulting problem to
obtain approximate solutions $(w_{\lambda,n},\vphi_{\lambda,n})$.
A priori estimates similar to ones performed below 
provide the compactness framework to pass to the limit as the 
Galerkin parameter $n\to \infty$ and as the regularization parameter $\lambda\to 0$, yielding a weak solution to the original
problem. Details of this procedure for constant mobility can be found
in \cite{CGM2024} (see also \cite[Theorem 4.8]{MB2019}).
In particular, the boundedness of $\vphi$ is proven using a classical argument from \cite{DebDet} (see also \cite[Theorem 4.8]{MB2019}) which yields
\begin{equation}\label{phi Linfty}
	\vphi \in L^{\infty}(\Omega \times (0,\infty)) \quad \text{and} \quad |\vphi(x,t)|<1 \ \text{a.e. in} \ \Omega, \ \forall \, t \in  (0,\infty).
\end{equation}
\textbf{A priori estimates}.
Testing \eqref{w = phi + ve mu} by $\vphi$ yields
\begin{equation*}
	(w, \vphi) = \ve \norm{\nabla \vphi}_{L^2(\Omega)}^2 + \ve (\Psi '(\vphi), \vphi) + \norm{\vphi}_{L^2(\Omega)}^2.
\end{equation*}
Exploiting the inequality
\begin{equation*}
	(\Psi '(\vphi), \vphi)\ge - \theta_0 \norm{\vphi}_{L^2(\Omega)}^2,
\end{equation*}
we get
\begin{equation*}
	\ve (\Psi '(\vphi), \vphi) + \norm{\vphi}_{L^2(\Omega)}^2 \ge (1- \ve \theta_0) \norm{\vphi}_{L^2(\Omega)}^2= 2c_1 \norm{\vphi}_{L^2(\Omega)}^2,
\end{equation*}
where $c_1 := \frac{1-\ve \theta_0}{2}>0$ by assumption \eqref{epps}.
Moreover, an application of H{\"o}lder and Young's inequalities yields $$|(w, \vphi)| \le c_1 \norm{\vphi}_{L^2(\Omega)}^2 + \frac{1}{4 c_1} \norm{w}_{L^2(\Omega)}^2.$$
Hence,
\begin{equation}\label{phi H1 w L2}
	\ve \norm{\nabla \vphi}_{L^2(\Omega)}^2+ c_1 \norm{\vphi}_{L^2(\Omega)}^2 \le \frac{1}{4 c_1} \norm{w}_{L^2(\Omega)}^2.
\end{equation}
Choosing $v=w$ in \eqref{w-weak}, we obtain
\begin{equation*}
	\frac{\ve}{2}\ddt \norm{w}_{L^2(\Omega)}^2 +(b(\vphi) \nabla w, \nabla w) = (b(\vphi) \nabla \vphi, \nabla w),
\end{equation*}
so, by H{\"o}lder and Young's inequalities we have
\begin{equation*}
	\frac{\ve}{2}\ddt \norm{w}_{L^2(\Omega)}^2 + \norm{\sqrt{b(\vphi)} \nabla w}_{L^2(\Omega)}^2 \le  \frac{1}{2} \norm{\sqrt{b(\vphi)} \nabla w}_{L^2(\Omega)}^2 + \frac{b_M}{2} \norm{\nabla \vphi}_{L^2(\Omega)}^2.
\end{equation*}
Therefore, exploiting \eqref{phi H1 w L2}, we infer that
\begin{equation*}
	\ve \ddt \norm{w}_{L^2(\Omega)}^2 + \norm{\sqrt{b(\vphi)} \nabla w}_{L^2(\Omega)}^2 \le \frac{b_M}{4c_1 \ve} \norm{w}_{L^2(\Omega)}^2.
\end{equation*}
An application of the standard Gr{\"o}nwall's Lemma on $[0,T]$, together with $w_0 \in L^2(\Omega)$, yields
\begin{equation}\label{first w regularity}
	w \in L^{\infty}(0, T; L^2(\Omega)), \quad \ \sqrt{b(\vphi)} \nabla w \in L^2(0, T; L^2(\Omega)).
\end{equation}
We thus learn from \eqref{phi H1 w L2} that 
\begin{equation}
	\label{phiH1}
\vphi\in L^{\infty}(0, T; H^1(\Omega)).
\end{equation} 
Since the mobility $b$ is bounded,  by comparison in the first equation, we also find
\begin{equation}
	\label{dabliu}
	\pt w\in  L^2(0,T; H^{-1}_{(0)}(\Omega)).
\end{equation}
Moreover, being the mobility non-degenerate, we have $b_m\|\nabla w\|^2_{L^2(\Omega)}\leq \|\sqrt{b(\vphi)} \nabla w\|^2_{L^2(\Omega)}$ and we also deduce that
$\nabla w \in L^2(0, T; L^2(\Omega))$, which, combined with
\eqref{first w regularity}, yields
\begin{equation}
	\label{w-H1}
	w \in L^2(0, T; H^1_{(m)}(\Omega)).
\end{equation}
Multiplying \eqref{w = phi + ve mu} by $-\Delta \vphi$ and exploiting the boundary conditions, we get
\begin{equation*}
	-(w, \Delta \vphi)= \ve \norm{\Delta \vphi}_{L^2(\Omega)}^2 + \ve(\Psi '' (\vphi) \nabla \vphi, \nabla \vphi) + \norm{\nabla \vphi}_{L^2(\Omega)}^2.
\end{equation*}
Since
\begin{equation*}
	(\Psi '' (\vphi) \nabla \vphi, \nabla \vphi) \ge -\theta_0 \norm{\nabla \vphi}_{L^2(\Omega)}^2
\end{equation*}
and
\begin{equation*}
	|(w, \Delta \vphi)| \le \frac{\ve}{2} \norm{\Delta \vphi}_{L^2(\Omega)}^2 + \frac{1}{2 \ve} \norm{w}_{L^2(\Omega)}^2,
\end{equation*}
it holds
\begin{equation}\label{phi H2 w L2 approx}
	\ve \norm{\Delta \vphi}_{L^2(\Omega)}^2 + 4c_1 \norm{\nabla \vphi}_{L^2(\Omega)}^2 \le \frac{1}{\ve} \norm{w}_{L^2(\Omega)}^2.
\end{equation}
In conclusion, recalling \eqref{first w regularity}, we have proved that $\vphi \in L^{\infty}(0, T; H^2(\Omega))$.
As a consequence, by comparison in \eqref{w = phi + ve mu}, we deduce that
\begin{equation}
	\label{F'-L2}
	F ' (\vphi)\in L^\infty(0, T; L^2(\Omega)),
\end{equation}
which, in turn, thanks to an argument in \cite{DebDet}, entails that
\begin{equation}
\label{phi Linfty-v2}
 |\vphi(x,t)|<1 \quad \text{a.e. in} \ \Omega,\  \forall \, t >0.
\end{equation}
Differentiating \eqref{w = phi + ve mu} with respect to time and multiplying by $\pt \vphi$ yields
\begin{equation} \label{diff eq 2}
	\l \pt w, \pt \vphi \r = \ve \norm{\nabla \pt \vphi}_{L^2(\Omega)}^2 + \ve( \Psi ''(\vphi) \pt \vphi, \pt \vphi) + \norm{\pt \vphi}_{L^2(\Omega)}^2 \ge \ve \norm{\nabla \pt \vphi}_{L^2(\Omega)}^2 + 2c_1 \norm{\pt \vphi}_{L^2(\Omega)}^2.
\end{equation}
Recalling that $\overline{\pt w}=0$, we have
\begin{equation}\label{eq 1 test pt w}
	\l \pt w, \pt \vphi \r \le \norm{\pt w}_{H^{-1}_{(0)}(\Omega)} \norm{\nabla \pt \vphi}_{L^2(\Omega)} \le \frac{1}{2 \ve} \norm{\pt w}_{H^{-1}_{(0)}(\Omega)}^2 + \frac{\ve}{2} \norm{\nabla \pt \vphi}_{L^2(\Omega)}^2.
\end{equation}
Therefore, we obtain
\begin{equation}\label{pt vhi H1 abd pt w dual}
	{\ve} \norm{\nabla \pt \vphi}_{L^2(\Omega)}^2 + 4c_1 \norm{\pt \vphi}_{L^2(\Omega)}^2 \le \frac{1}{\ve} \norm{\pt w }^2_{H^{-1}_{(0)}(\Omega)},
\end{equation}
which implies $\pt \vphi \in L^2(0,T ; H^1(\Omega))$.

\medskip
\noindent In order to complete the proof of part (A) of Theorem \ref{ws thm} we now pursue uniform in time estimates, for which we turn to the equivalent formulation \eqref{CH1-p}-\eqref{CH2-p}, having set
$$\mu = \frac{1}{\ve} (w- \vphi).$$
\textbf{Uniform in time estimates}.
First of all, taking $v = \mu$ as a test function in \eqref{weak formulation perturbed Cahn--Hilliard equation}, by using \eqref{mu perturbed Cahn--Hilliard equation}, and exploiting the standard chain rule in $L^2(0, T; H^1(\Omega)) \cap H^1(0, T; (H^1(\Omega))')$, for every $T>0$ we get
\begin{equation}\label{"Energy Balance"}
	\ddt\left( \frac12 \norm{\nabla \vphi}^2_{L^2(\Omega)}+ \frac{\ve}{2} \norm{\mu}_{L^2(\Omega)}^2 + \int_{\Omega} \Psi(\vphi) \, \d x \right) + \int_{\Omega} b(\vphi)|\nabla \mu |^2 \, \d x =0
\end{equation}
almost everywhere in $(0,\infty)$.
Thus, recalling that $\vphi_0 \in H^2(\Omega)$, with $\norm{\vphi_0}_{L^\infty(\Omega)}\leq 1$ and $\mu_0 \in L^2(\Omega)$, we deduce
\begin{align}
	\label{nabla phi L infty L2 uniform bound}
	&\norm{\nabla \vphi}_{L^\infty(0,\infty;L^2(\Omega))} \le C,\\
	\label{mu L infty L2 uniform bounds}
	&\sqrt{\ve} \norm{\mu}_{L^\infty(0,\infty;L^2(\Omega))} \le C
\end{align}
and
\begin{equation}\label{nabla mu L2 L2 uniform bound}
	\int_0^\infty \norm{\sqrt{b(\vphi)} \nabla \mu}_{L^2(\Omega)}^2 \, \d s \le C,
\end{equation}
where here and throughout this proof $C$ stands for a positive constant depending on the norm of the initial data but independent of time and the parameter $\ve$.
By the non-degeneracy of $b$ we also get
\begin{equation}
	\label{nabla mu L2L2 uniform in time}
	\norm{\nabla \mu}_{L^2(0, \infty; L^2(\Omega))}^2 
	\le \frac{1}{b_m} \int_0^\infty \norm{\sqrt{b(\vphi)} \nabla \mu}^2_{L^2(\Omega)}  \, \d s \le C.
\end{equation}
Next, exploiting \eqref{nabla phi L infty L2 uniform bound} and the constraint $\norm{\vphi(t)}_{L^\infty(\Omega)} \le 1$, it follows
\begin{equation}
	\label{phi L infty H1 uniform bound}
	\norm{\vphi}_{L^\infty(0,\infty;H^1(\Omega))} \le C,
\end{equation}
as well as
\begin{equation}
	\label{w L infty L2 uniform bound}
	\norm{w}_{L^\infty(0,\infty;L^2(\Omega))} \le C,
\end{equation}
on account of \eqref{mu L infty L2 uniform bounds}, since $w=\vphi+\ve \mu$. Using this in \eqref{phi H2 w L2 approx} we immediately learn that
  \begin{equation}
 	\label{phi L infty H2 unif bound}
 	\ve\norm{\Delta\vphi}_{L^\infty(0,\infty;L^2(\Omega))} \le C,
 \end{equation}
which, by comparison in \eqref{CH2-p}, also implies $\ve\norm{F'(\vphi)}_{L^\infty(0, \infty; L^2(\Omega))} \le C$. Besides, collecting
\eqref{nabla mu L2L2 uniform in time}, \eqref{phi L infty H1 uniform bound} and \eqref{w L infty L2 uniform bound},  we have
\begin{align}
	\label{w L2 uloc H1}
	\int_t^{t+1} \norm{w}_{H^1(\Omega)}^2 \d s &= \int_t^{t+1} \norm{w}_{L^2(\Omega)}^2 \d s+ \int_t^{t+1} \norm{\nabla w}_{L^2(\Omega)}^2 \d s\nonumber \\
	&\le \norm{w}_{L^\infty(0, \infty; L^2(\Omega))}^2+ 2 \ve^2 \norm{\nabla \mu}_{L^2(0, \infty; L^2(\Omega))}^2+ 2\norm{\nabla \vphi}_{L^\infty(0, \infty; L^2(\Omega))}^2 \le C,
\end{align}
yielding $w\in L^2_{\rm uloc}([0, \infty); H^1_{(m)}(\Omega))$.
Finally, owing to \eqref{nabla mu L2 L2 uniform bound}, a comparison in \eqref{CH1-p} yields
\begin{equation}
	\label{pt w L2 dual uniform bound}
	\int_0^\infty \norm{\pt\vphi+\ve\pt\mu}_{H^{-1}_{(0)}(\Omega)}^2 \, \d s \le C,
\end{equation}
namely  $\pt w\in L^2(0, \infty; H^{-1}_{(0)}(\Omega))$.
As a consequence, reasoning as in \eqref{pt vhi H1 abd pt w dual}, we infer $ \pt \vphi \in L^2(0, \infty; H^1(\Omega))$, with
 \begin{equation}
 	\label{pt vphi L2 H1}
 	\ve \norm{\pt \vphi}_{L^2(0, \infty; H^1(\Omega))}\leq C.
 \end{equation}
Collecting all the previous estimates we then accomplish statement (A).
\begin{remark}
\label{lambda independence 2}
Notice that we exploited the non-degeneracy of the mobility only in the proof of \eqref{w-H1}, \eqref{nabla mu L2L2 uniform in time} and \eqref{w L2 uloc H1}.
Hence, all the other estimates hold true even when no positive lower bound for the mobility is available.
Furthermore, it is apparent from the proof that some of the bounds we obtained are uniform also with respect to the parameter $\ve$.
We will exploit and deepen this fact in Section \ref{S-convergenza}.
\end{remark}

\noindent \textbf{Uniqueness of weak solution}.
Let $(w_i, \vphi_i)$, $i=1,2$, be two weak solutions to \eqref{w-weak}--\eqref{w = phi + ve mu} with initial conditions $w_1^0, w_2^0$ respectively, such that $\overline{w_1^0}=\overline{w_2^0}$. 
We will denote by  $C$ a positive constant depending on $\ve$, the structural parameters,
$\norm{w_1^0}_{L^2(\Omega)}, \norm{w_2^0}_{L^2(\Omega)}$ and the mobility $b$, but independent of time.

Define $W=w_1 - w_2$ and $\Phi=\vphi_1-\vphi_2$. We observe that they satisfy
\begin{equation}\label{w diff 1}
W = -\ve \left( \G_{\vphi_1} \pt w_1 - \G_{\vphi_2} \pt w_2 \right) + \Phi - \overline{\Phi} = -\ve \left( \G_{\vphi_1} \pt W + \left( \G_{\vphi_1} - \G_{\vphi_2}\right) \pt w_2 \right) + \Phi -\overline{\Phi}
\end{equation}
and
\begin{equation}\label{w diff 2}
W= \ve( - \Delta \Phi + \Psi'(\vphi_1) - \Psi'(\vphi_2)) + \Phi.
\end{equation}
Testing \eqref{w diff 2} by $ \Phi$, we have
$$(W, \Phi) = \ve \norm{\nabla \Phi}^2_{L^2(\Omega)}+ \ve( \Psi'(\vphi_1) - \Psi'(\vphi_2), \Phi) + \norm{\Phi}^2_{L^2(\Omega)}.$$
The left-hand side is estimated as
$$|(W, \Phi)| \le  \norm{W}_{H^{-1}_{(0)}(\Omega)} \norm{\nabla \Phi}_{L^2(\Omega)}$$
while
$$ \ve( \Psi'(\vphi_1) - \Psi'(\vphi_2), \Phi) + \norm{\Phi}^2_{L^2(\Omega)} = (1- \ve \theta_0) \norm{\Phi}^2_{L^2(\Omega)} + \ve (F'(\vphi_1) - F'(\vphi_2), \Phi) \ge 2 c_1 \norm{\Phi}^2_{L^2(\Omega)},$$
where $c_1=\frac{1- \ve \theta_0}2>0$. Then, an application of Young's inequality yields
\begin{equation}\label{norm phi e norm w}
\norm{\Phi}^2_{H^1(\Omega)} \le C \norm{W}_{H^{-1}_{(0)}(\Omega)}^2.
\end{equation}
Testing \eqref{w diff 1} with $W$, we find
\begin{equation*}
\ve  \left( \G_{\vphi_1} \pt W, W \right)+\norm{W}_{L^2(\Omega)}^2 = - \ve\left(  \left( \G_{\vphi_1} - \G_{\vphi_2}\right) \pt w_2, W \right) + (\Phi , W),
\end{equation*}
and, on account of \eqref{pt vphi L2 H1}, the first term reads as (cfr. \cite[Eqn. (2.22)]{BB1999}, \cite[Eqn. (3.38)]{CGGG2025})
\begin{equation}
	\label{integration by part formula}
	2\left( \G_{\vphi_1} \pt W, W  \right) = \ddt \norm{\sqrt{b(\vphi_1) }\nabla \G_{\vphi_1} W}_{L^2(\Omega)}^2 + \left( b'(\vphi_1) \pt \vphi_1, |\nabla \G_{\vphi_1} W|^2 \right).
\end{equation}
Therefore, we infer
\begin{equation}\label{diff inq 1}
 \frac{\ve}{2} \ddt \norm{\sqrt{b(\vphi_1) }\nabla \G_{\vphi_1} W}_{L^2(\Omega)}^2+\norm{W}_{L^2(\Omega)}^2= (\Phi, W)+ I_1+I_2,
\end{equation}
where
$$ 
I_1 = - \frac{\ve}{2} \left( b'(\vphi_1) \pt \vphi_1, |\nabla \G_{\vphi_1} W|^2 \right)\quad\text{and}\quad
I_2= - \ve \left(  \left( \G_{\vphi_1} - \G_{\vphi_2}\right) \pt w_2, W \right).
$$
By the previous computations and \eqref{equivalence of norm},
$$|(\Phi, W)|\leq C\norm{W}^2_{H^{-1}_{(0)}(\Omega)}\leq C\norm{\sqrt{b(\vphi_1) }\nabla \G_{\vphi_1} W}_{L^2(\Omega)}^2.
$$
Exploiting the boundedness of $b$ and the interpolation inequality $\norm{v}_{L^{\frac{12}{5}}(\Omega)} \le C \norm{v}_{L^2(\Omega)}^{\frac{3}{4}}\norm{v}_{H^1(\Omega)}^{\frac{1}{4}}$,
we obtain
$$|I_1| \le C \norm{ \pt \vphi_1}_{L^6(\Omega)} \norm{\nabla \G_{\vphi_1} W}^2_{L^\frac{12}{5}(\Omega)}\le C \norm{ \pt \vphi_1}_{H^1(\Omega)} \norm{\nabla \G_{\vphi_1} W}^{\frac32}_{L^2(\Omega)}\norm{\G_{\vphi_1} W}^{\frac12}_{H^2(\Omega)}.$$
Recalling \eqref{H2 norm of G_q-3d}, in light of the uniform control $\ve\vphi_1\in L^\infty(0,\infty; H^2(\Omega))$, we have
\begin{equation}
\label{H2 norm of G_qbis}
	\norm{\G_{\vphi_1} W}_{H^2(\Omega)} \le C \left( \norm{W}_{L^2(\Omega)}+ \norm{{\nabla\G_{\vphi_1} W}}_{L^2(\Omega) } \right).
\end{equation}
Therefore, we deduce that
\begin{align}
\label{control I1}|I_1| & \le
C \norm{ \pt \vphi_1}_{H^1(\Omega)} \norm{\nabla \G_{\vphi_1} W}^{\frac32}_{L^2(\Omega)}\norm{W}^{\frac12}_{L^2(\Omega)}
+C \norm{ \pt \vphi_1}_{H^1(\Omega)} \norm{\nabla \G_{\vphi_1} W}^{2}_{L^2(\Omega)}\nonumber\\
&  \leq \frac14 \norm{W}^2_{L^2(\Omega)}+ C\norm{ \pt \vphi_1}^{\frac43}_{H^1(\Omega)} \norm{\nabla \G_{\vphi_1} W}^2_{L^2(\Omega) }
+C \norm{ \pt \vphi_1}_{H^1(\Omega)} \norm{\nabla \G_{\vphi_1} W}^{2}_{L^2(\Omega)}  \nonumber\\
& \leq \frac14 \norm{W}^2_{L^2(\Omega)}+ C( 1+\norm{ \pt \vphi_1}^2_{H^1(\Omega)} )
\norm{\sqrt{b(\vphi_1)}\nabla \G_{\vphi_1} W}^2_{L^2(\Omega) }  .
\end{align}
The estimate of $I_2$ relies on the definition of $\G_q$ and its symmetry.  Indeed,
\begin{equation}\label{I3 est}
\begin{split}
|I_2| &= | \ve \left(  \left( \G_{\vphi_1} - \G_{\vphi_2}\right) \pt w_2, W \right)| = |\l \left( \G_{\vphi_1} - \G_{\vphi_2}\right) W,\ve\pt w_2 \r| \\
& = |\left( b(\vphi_2) \nabla \G_{\vphi_2} (\ve \pt w_2), \nabla \left(\G_{\vphi_1} - \G_{\vphi_2}\right) W\right)|\\
&= |\left(\nabla \G_{\vphi_2} (\ve \pt w_2), (b(\vphi_2)-b(\vphi_1)) \nabla \G_{\vphi_1} W \right)|\\
& \le \norm{\nabla \G_{\vphi_2} (\ve \pt w_2)}_{L^2(\Omega)} \norm{ (b(\vphi_2)-b(\vphi_1)) \nabla \G_{\vphi_1} W}_{L^2(\Omega)}.
\end{split}
\end{equation}
Since $\nabla \G_{\vphi_2} (\ve \pt w_2 )= \nabla( \vphi_2-w_2)=-\ve\nabla \mu_2$, the first norm is bounded as
$$\norm{\nabla \G_{\vphi_2} (\ve \pt w_2)}_{L^2(\Omega)} \le C\norm{\nabla \mu_2}_{L^2(\Omega)}.$$
Concerning the second term in \eqref{I3 est}, the interpolation inequality $\|v\|_{L^3(\Omega)}^2\le C\|v\|_{L^2(\Omega)}\|v\|_{H^1(\Omega)}$ yields
\begin{equation*}
		\norm{(b(\vphi_2)-b(\vphi_1))\nabla \G_{\vphi_1} W}_{L^2(\Omega)} \le C \norm{\Phi}_{H^1(\Omega)} \norm{\nabla \G_{\vphi_1}W}_{L^2(\Omega)}^\frac12 \norm{\G_{\vphi_1}W}_{H^2(\Omega)}^\frac12.
\end{equation*}
Therefore, recalling \eqref{equivalence of norm}, \eqref{norm phi e norm w} and \eqref{H2 norm of G_qbis}, we find
\begin{align*}
|I_2| &\le C {\norm{\nabla \mu_2}_{L^2(\Omega)}} \norm{\Phi}_{H^1(\Omega)} \norm{\nabla\G_{\vphi_1}W}_{L^2(\Omega)}^\frac12 \norm{\G_{\vphi_1}W}_{H^2(\Omega)}^\frac12\\
&\le C{\norm{\nabla \mu_2}_{L^2(\Omega)}} \norm{\nabla\G_{\vphi_1}W}_{L^2(\Omega)}^\frac32\norm{\G_{\vphi_1}W}_{H^2(\Omega)}^\frac12\\
& \leq \frac14 \norm{W}^2_{L^2(\Omega)}+ C  \norm{\nabla \mu_2}^{\frac43}_{L^2(\Omega)} \norm{{\nabla\G_{\vphi_1} W}}^2_{L^2(\Omega) }
+ C  \norm{\nabla \mu_2}_{L^2(\Omega)} \norm{{\nabla\G_{\vphi_1} W}}^2_{L^2(\Omega) }\\
& \leq \frac14 \norm{W}^2_{L^2(\Omega)} + C (1 +
\norm{\nabla \mu_2}^2_{L^2(\Omega)}) \norm{\sqrt{b(\vphi_1) }\nabla\G_{\vphi_1} W}^2_{L^2(\Omega)}.
\end{align*}
Therefore, \eqref{diff inq 1} yields the final differential inequality
\begin{align*}
\ve \ddt\norm{\sqrt{b(\vphi_1) }\nabla \G_{\vphi_1} W}_{L^2(\Omega)}^2
		& \le h(t)\norm{\sqrt{b(\vphi_1) }\nabla \G_{\vphi_1} W}_{L^2(\Omega)}^2 ,
\end{align*}
where
$$
h(t) = C(1+\norm{ \pt \vphi_1}_{H^1(\Omega)}^2 + \norm{\nabla \mu_2}_{L^2(\Omega)}^2)\in L^1_{\rm{uloc}}([0, \infty))
$$
in view of  \eqref{nabla mu L2L2 uniform in time} and \eqref{pt vphi L2 H1}.
Applying Gr{\"o}nwall's inequality, recalling the equivalence \eqref{equivalence of norm} and the control \eqref{norm phi e norm w},
we find the continuous dependence estimate
$$
\norm{W(t)}^2_{H^{-1}_{(0)}(\Omega)}+\norm{\Phi(t)}^2_{H^1(\Omega)}\leq Ce^{Ct}\norm{W(0)}^2_{H^{-1}_{(0)}(\Omega)} ,\quad \forall \, t\geq 0.
$$
Finally, since for equal initial data $w_1^0= w_2^0$ we have $W(0)=0$, uniqueness of the weak solution follows by
$\norm{W(t)}_{H^{-1}_{(0)}(\Omega)}^2=0$ for any $t\geq 0$.

\begin{remark}\label{rem-data}
In Theorem \ref{ws thm} we assumed that $w_0 \in L^2(\Omega)$, which implies that $\vphi_0 \in H^2(\Omega)$. So, we have a stronger initial datum compared with the classical Cahn--Hilliard equation, where $\vphi_0 \in H^1(\Omega)$. The existence of weak solutions for initial data $\vphi_0\in H^1(\Omega)$ remains an open question, even in the case of constant mobility \eqref{m-const}.
\end{remark}

\noindent
\textbf{Further regularity and separation from the pure phases.} In order to prove that
$w\in L^\infty(\Omega\times (\tau,\infty))$, we shall exploit the Nash--De Giorgi theory for parabolic equation as in \cite{DiBenedettoDegenerate}.
To this aim, we recall the following well-known result.
\begin{lemma}\label{De Giorgi Lemma}
	Let $\{Z_n\}_{n \in \mathbb{N}} \subset \R_+$ satisfy the relation
	\begin{equation*}
		Z_{n+1}\le \widetilde{C}\beta^n Z_n^{1+ \alpha},
	\end{equation*}
	for some $\widetilde{C}>0,\ \beta>1$ and $\alpha>0$. Assume that $Z_0 \le \widetilde{C}^{-\frac{1}{\alpha}} \beta^{- \frac{1}{\alpha^2}}$. Then we have
	\begin{equation*}
		Z_n \le Z_0 \beta^{-\frac{n}{\alpha}}, \quad \forall \, n \ge 1.
	\end{equation*}
	In particular, $Z_n \to 0$ as $n \to \infty$.
\end{lemma}
Now, for any $\tau >0$ consider the time step $2 \tau$. For all $n \in \mathbb{N}$,  $0 \le t_\ast<\infty$ we define the sequence of increasing times and rising levels
	\begin{align*}
		t_n = \tau \left(1- \frac{1}{2^n}\right)+ t_\ast,\qquad
		M_n = L \left(1- \frac{1}{2^n}\right),
	\end{align*}
	where $L>1$ will be chosen later. We consider
	\begin{equation*}
		y_n = (w-M_n)_+ \quad \text{and} \quad x_n = (w+M_n)_-,
	\end{equation*}
where $(f)_+:=\sup\{0,f\}$ and $(f)_-:=\sup\{0,-f\}$.
	Our goal is to apply Lemma \ref{De Giorgi Lemma} to the sequences $Y_n$ and $X_n$ defined as 	
	\begin{align*}
		Y_n =  \int_{t_n}^{2 \tau + t_\ast} \int_{\Omega}y_n^2 \,\d x \, \d s \quad \text{and} \quad
		X_n = \int_{t_n}^{2 \tau + t_\ast} \int_{\Omega}x_n^2 \,\d x \, \d s,
	\end{align*}
which will provide $Y_n \to 0$ and $X_n \to 0$ in the limit $n\to\infty$. Therefore, observing that
$$Y_n \to \int_{\tau+ t_\ast}^{2 \tau + t_\ast}\int_{\Omega}(w-L)_+^2 \,\d x \, \d s,  \quad X_n \to \int_{\tau+t_\ast}^{2 \tau + t_\ast}\int_{\Omega}(w+L)_-^2\, \d x \, \d s \quad \text{as}\ n \to \infty,$$
this will yield $w \in L^{\infty}(\Omega \times [\tau+t_\ast, 2 \tau +t_\ast]), \ \forall \, t_\ast \ge 0$. Since it turns out that this bound does not depend on $t_\ast$, iterating the argument we will obtain a global in time bound for $w \in L^{\infty}(\Omega \times [\tau, \infty))$, for any $\tau >0$. The proof is divided into two steps.

\smallskip
\noindent \textbf{Step 1: De Giorgi's scheme.} From this point onward $C$ will denote a positive constant that depends on the parameters of the system, $\tau$ and $\ve$. We begin by defining the following family of cut-off functions
\begin{equation*}
	\eta_n (s) = \begin{cases}
		1 \quad & t_{n+1} \le s \le 2 \tau + t_\ast, \\
		\frac{s-t_n}{t_{n+1}-t_n} \ & t_n < s < t_{n+1},\\
		0 & 0 \le s \le t_n.
			\end{cases}
\end{equation*}
For all $z >0$, we also define
\begin{equation}
\label{A_n set}
A_n(z) := \{x \in \Omega :\ w(x, z)>M_n\}
\end{equation}
and observe that $\chi_{A_n}= \chi_{\lbrace y_{n}>0\rbrace }$ where, $\forall A \subset \R^d$, $\chi_A$ is the characteristic function of the set $A$, namely
\begin{equation*}
	\chi_A(x)= \begin{cases}
		1 \quad \text{ if } x \in A,\\
		0 \quad \text{ otherwise}.
	\end{cases}
\end{equation*}
By using \eqref{w-weak} with $v=y_{n+1} \eta_n \in H^1(\Omega)$, we get
\begin{equation}
\label{DG-test}
\langle \ve \pt w, y_{n+1} \rangle \eta_n+ (b(\vphi)\nabla w, \nabla y_{n+1})\eta_n = (b(\vphi) \nabla \vphi, \nabla y_{n+1})\eta_n.
	\end{equation}
Since $w\in L^2(0,T;H^1(\Omega))\cap H^1(0,T;(H^1(\Omega))')$, 
by the chain rule applied to $\l\pt w,(w-M_{n+1})_+\r$, we observe that
\begin{equation*}
\langle \ve \pt w, y_{n+1} \rangle \eta_n = \frac{\ve}{2} \ddt \int_{\Omega} y_{n+1}^2 \eta_n \,\d x - \frac{\ve}{2} \int_{\Omega}y_{n+1}^2 \pt \eta_n \,\d x
	\end{equation*}
almost everywhere on $(0,\infty)$.
We also notice that
	\begin{equation*}
		(b(\vphi)\nabla w, \nabla y_{n+1})\eta_n = \int_{\Omega} b(\vphi)|\nabla y_{n+1}|^2 \eta_n \,\d x.
	\end{equation*}
Moreover, we estimate the right-hand side of \eqref{DG-test} as follows
\begin{align*}
(b(\vphi) \nabla \vphi, \nabla y_{n+1})\eta_n 
&= \int_{\Omega} b(\vphi) \nabla \vphi \cdot \nabla y_{n+1}\eta_n \,\d x 
 = \int_{\Omega} b(\vphi) \chi_{A_{n+1}} \nabla\vphi\cdot \nabla y_{n+1} \eta_n \,\d x\nonumber \\
		& \le\frac12 \int_\Omega b(\vphi) |\nabla \vphi|^2 \chi_{A_{n+1}} \eta_n \,\d x
		+ \frac12\int_{\Omega}b(\vphi) |\nabla y_{n+1}|^2 \eta_n \,\d x.	
	\end{align*}
	Therefore, we infer
	\begin{align*}
		\ve\ddt \int_{\Omega}y_{n+1}^2 \eta_n \,\d x + \int_{\Omega}b(\vphi)|\nabla y_{n+1}|^2 \eta_n \,\d x &\le \int_{\Omega}b(\vphi) |\nabla\vphi|^2 \chi_{A_{n+1}} \eta_n \,\d x +\ve	\int_{\Omega}y_{n+1}^2 \pt \eta_n \,\d x.
	\end{align*}
Next, integrating between 0 and any $z\in [t_{n+1}, 2 \tau + t_\ast]$, and exploiting the definition of $\eta_n$, we get
	\begin{align*}
		&\ve \int_\Omega y_{n+1}^2 \,\d x + \int_{t_{n+1}}^{z} \int_\Omega b(\vphi) |\nabla y_{n+1}|^2  \,\d x \d s \nonumber \\
		&\quad \le \int_0^{z} \int_{\Omega}b(\vphi)|\nabla \vphi|^2 \chi_{A_{n+1}} \eta_n \,\d x \, \d s + \ve \int_0^{z} \int_{\Omega}y_{n+1}^2 \pt \eta_n \,\d x \, \d s \nonumber \\
		&\quad \le \int_0^{2 \tau + t_\ast} \int_{\Omega}b(\vphi)|\nabla \vphi|^2 \chi_{A_{n+1}} \eta_n \,\d x \, \d s + \ve \int_0^{2 \tau + t_\ast} \int_{\Omega}y_{n+1}^2 \pt \eta_n \,\d x\, \d s.
	\end{align*}
	Taking the supremum on $z$ in the time interval $[t_{n+1}, 2 \tau + t_\ast]$, we thus have
	\begin{align}
		\label{first de giorgi inequality}
		&\ve\sup_{t_{n+1}\le z \le 2 \tau + t_\ast} \int_{\Omega}y_{n+1}^2 \,\d x  + \int_{t_{n+1}}^{2 \tau + t_\ast} \int_{\Omega} b(\vphi) |\nabla y_{n+1}|^2  \,\d x \,\d s \nonumber\\
		&\quad \le \int_0^{2 \tau + t_\ast} \int_{\Omega}b(\vphi)|\nabla \vphi|^2 \chi_{A_{n+1}} \eta_n \,\d x \,\d s + \ve \int_0^{2 \tau + t_\ast} \int_{\Omega}y_{n+1}^2 \pt \eta_n \,\d x \, \d s.
	\end{align}
	Now we estimate the right-hand side. Applying the definition of $\eta_n$, \eqref{A_n set} and H{\"o}lder's inequality, we obtain
	\begin{align*}
		&\int_0^{2 \tau + t_\ast} \int_{\Omega}b(\vphi) |\nabla \vphi|^2\chi_{A_{n+1}} \eta_n\,\d x \, \d s \nonumber \\
		&\quad \le \int_{t_n}^{2 \tau + t_\ast} \int_{\Omega}b(\vphi) |\nabla \vphi|^2\chi_{A_{n+1}} \,\d x \, \d s  \nonumber \\
		&\quad\le \left( \int_{t_n}^{2 \tau + t_\ast} \int_{\Omega}\big(b(\vphi)|\nabla\vphi|^2\big)^3 \,\d x \,\d s \right)^\frac13 \left(\int_{t_n}^{2 \tau + t_\ast} \int_{\Omega}\chi_{A_{n+1}}^\frac32 \,\d x \,\d s \right)^\frac23 \nonumber\\
		&\quad\le C \left(2\tau-\tau\left(1- \frac{1}{2^n}\right)\right)^\frac13 \norm{\vphi}_{L^{\infty}(0, \infty; H^2(\Omega))}^2\left(\int_{t_n}^{2 \tau + t_\ast}\int_{\Omega}
		\left|\chi_{A_{n+1}}\right|^\frac32\,\d x\, \d s\right)^\frac23 \nonumber\\
		&\quad \le C \left(\tau \left(1+ \frac{1}{2^n}\right)\right)^\frac13 \norm{\vphi}_{L^{\infty}(0, \infty; H^2(\Omega))}^2 \left(\int_{t_n}^{2 \tau + t_\ast} |A_{n+1}(s)|\,\d s\right)^\frac23\!\!.
	\end{align*}
Then we control $|A_{n+1}(z)|$ as in \cite{DiBenedettoDegenerate}, writing
	\begin{align}
		\label{A_n estimate}
		\frac{L^2}{2^{2(n+1)}}\int_{t_n}^{2 \tau + t_\ast}|A_{n+1}(s)| \,\d s &= (M_{n+1}-M_n)^2\int_{t_n}^{2 \tau + t_\ast}|A_{n+1}(s)| \,\d s \nonumber \\
		&\le \int_{t_n}^{2 \tau + t_\ast} \int_{\Omega} y_n^2 \chi_{A_{n+1}} \,\d x \, \d s 
\le \int_{t_n}^{2 \tau + t_\ast} \int_{\Omega}y_n^2 \,\d x \, \d s.
	\end{align}
	Moreover, recalling the definition of $\eta_n$ we have
	\begin{equation*}
		\pt \eta_n\le \frac{1}{t_{n+1}-t_n} =\frac{2^{n+1}}{\tau}, \quad \forall \, n \in \mathbb{N}.
	\end{equation*}
Thus, combining \eqref{first de giorgi inequality}, \eqref{A_n estimate}, the monotonicity of the sequence $\{M_n\}$ and the above inequality, we get
	\begin{align}
		\label{equation est for De Giorgi}
		&\sup_{t_{n+1} \le z \le 2 \tau + t_\ast} \ve \int_{\Omega} y_{n+1}^2 \,\d x + \int_{t_{n+1}}^{2 \tau + t_\ast} \int_{\Omega} b(\vphi)|\nabla y_{n+1}|^2 \,\d x \, \d s \nonumber \\
		&\quad \le C\frac{2^{\frac43 (n+1)}}{L^\frac43} \left(\int_{t_n}^{2 \tau + t_\ast} \int_{\Omega} y_n^2 \,\d x \d s\right)^\frac23 +\ve\frac{2^{n+1}}{\tau}\int_{t_n}^{2 \tau + t_\ast}\int_{\Omega} y_n^2\, \d x \, \d s \nonumber\\
		&\quad \le C\frac{2^{\frac43 (n+1)}}{L^\frac43} Y_n^\frac23 + \ve\frac{2^{n+1}}{\tau} Y_n.
	\end{align}

	We are now ready to apply the De Giorgi iterative scheme. We observe that, due to \eqref{ws regularity}, it holds
	\begin{equation}
		\label{Yn bound}
		Y_n \le C, \quad \forall n \in \mathbb{N},
	\end{equation}
	where $C$ is independent of $t_\ast$.
	Applying H{\"o}lder's inequality and \eqref{A_n estimate}, we find
	\begin{align}
		\label{Yn estimate}
		Y_{n+1} &\le \left(\int_{t_{n+1}}^{2 \tau + t_\ast} \int_{\Omega}y_{n+1}^\frac{10}{3} \,\d x \, \d s\right)^\frac35\left(\int_{t_{n+1}}^{2 \tau + t_\ast}|A_{n+1}(s)|\,\d s\right)^\frac25 \nonumber\\
		& \le \left(\int_{t_{n+1}}^{2 \tau + t_\ast} \int_{\Omega}y_{n+1}^\frac{10}{3} \,\d x \, \d s\right)^\frac35 \left(\frac{2^{2(n+1)}}{L^2} Y_n\right)^\frac25\!\!.
	\end{align}
Recalling the interpolation inequality
		$\norm{v}_{L^\frac{10}{3}(\Omega)} \le C \norm{v}_{L^2(\Omega)}^\frac25 \norm{v}_{H^1(\Omega)}^\frac35$
and \eqref{equation est for De Giorgi}, we infer
	\begin{align}
		\label{L 10 3 norm estimate}
		\int_{t_{n+1}}^{2 \tau + t_\ast} \int_{\Omega} y_{n+1}^\frac{10}{3} \,\d x \, \d s & \le C\int_{t_{n+1}}^{2 \tau + t_\ast} \norm{y_{n+1}}_{L^2(\Omega)}^\frac43 \norm{y_{n+1}}_{H^1(\Omega)}^2 \,\d s\nonumber \\
		& \le C \left(\sup \limits_{t_{n+1}\le z \le 2 \tau + t_\ast} \int_{\Omega} y_{n+1}^2\, \d x \right)^\frac53 
		\nonumber\\
		&\quad + C\left(\sup \limits_{t_{n+1}\le z \le 2 \tau + t_\ast} \int_{\Omega} y_{n+1}^2 \,\d x\right)^\frac23 \int_{t_{n+1}}^{2 \tau + t_\ast} \norm{\nabla y_{n+1}}_{L^2(\Omega)}^2 \,\d s    \nonumber \\
		& \le C \left(\frac{2^{\frac43 (n+1)}}{L^\frac43}Y_n^\frac23 + \frac{2^{n+1}}{\tau}Y_n \right)^\frac53\!\!,
	\end{align}
	where $C>0$ is again independent of $t_\ast$.
	Thus, combining \eqref{Yn bound}, \eqref{Yn estimate} and \eqref{L 10 3 norm estimate}, we arrive at
	\begin{align}
		\label{De Giorgi second estimate}
		Y_{n+1} &\le C \left(\frac{2^{\frac43 (n+1)}}{L^\frac43} Y_n^\frac23 + \frac{2^{n+1}}{\tau}Y_n\right)\left(\frac{2^{2(n+1)}}{L^2}Y_n\right)^\frac25 \nonumber \\
		 & \le C \left(\frac{2^{\frac43 (n+1)}}{L^\frac43}\frac{2^{\frac45 (n+1)}}{L^\frac45}Y_n^{1 + \frac{1}{15}}+ \frac{2^{n+1}}{\tau}\frac{2^{\frac45 (n+1)}}{L^\frac45} Y_n^{1+ \frac25}\right) \nonumber\\
		 &\le C \left(\frac{2^{\frac43 (n+1)}}{L^\frac43}\frac{2^{\frac45 (n+1)}}{L^\frac45}+ \frac{2^{n+1}}{\tau}\frac{2^{\frac45 (n+1)}}{L^\frac45}\right) Y_n^{1 + \frac{1}{15}} \nonumber\\
		 & \le C \frac{2^{\frac{32}{15}(n+1)}}{L^\frac45}\left(\frac{1}{\tau}+ \frac{1}{L^\frac43}\right)Y_n^{1+\frac{1}{15}}.
	\end{align}	
	Since $L \ge 1$, we end up with
	\begin{equation}
		\label{iterative estimate}
		Y_{n+1}\le C \ \frac{2^{\frac{32}{15}(n+1)}}{L^\frac45}\left(\frac{1}{\tau}+ 1\right)Y_n^{1+\frac{1}{15}}.
	\end{equation}

\medskip
\noindent \textbf{Step 2: Iteration argument.} In order to apply Lemma \ref{De Giorgi Lemma}, with $\widetilde{C}= C \frac{1}{L^\frac45}(1+\frac{1}{\tau})2^\frac{32}{15}$, $\beta= 2^\frac{32}{15}$ and $\alpha= \frac{1}{15}$ we need $L$ to satisfy
	\begin{equation*}
		Y_0 \leq  \int_{t_\ast}^{2 \tau + t_\ast} \int_{\Omega} |w|^2 \,\d x \, \d s \le C^{-15} 2^{-512}L^{12}\left(1+\frac{1}{\tau}\right)^{-15}\!\!.
	\end{equation*}
	Hence, since $w \in L^\infty(0, \infty; L^2(\Omega))$ by \eqref{ws regularity}, it suffices to choose $L>1$ such that
	\begin{equation*}
		\norm{w}_{L^\infty(0, \infty; L^2(\Omega))}^2 \le C^{-15}2^{-512}L^{12}\frac{1}{2 \tau}\left(1+\frac{1}{\tau}\right)^{-15}\!\!.
	\end{equation*}
Therefore, $L$ depends on $\tau$ and is independent of $t_\ast$. At this point, an application of Lemma \ref{De Giorgi Lemma} yields
	\begin{equation*}
		\lim\limits_{n \to \infty}Y_n= \int_{\tau+ t_\ast}^{2 \tau + t_\ast} \int_{\Omega}|(w- L)_+|^2 \,\d x \, \d s =0,
	\end{equation*}
	that is, that $\norm{(w-L)_+}_{L^2(\Omega \times (\tau + t_\ast, 2 \tau+ t_\ast))} =0$. Repeating the same argument for $X_n$ we also infer that $\norm{(w+L)_-}_{L^2(\Omega \times (\tau + t_\ast, 2 \tau+ t_\ast))} =0$.
Finally this implies the desired conclusion
$$ 
w \in L^{\infty}(\Omega\times[\tau+t_\ast, 2 \tau +t_\ast]), \quad \forall \, t_\ast \ge 0.
$$
	Therefore, choosing $t_\ast=0$ we infer
	$$w \in L^{\infty}(\Omega \times [\tau, 2 \tau]).$$
Now, since $L$ does not depend on $t_\ast$, we can choose $t_\ast= \tau $ to obtain $w \in L^{\infty}(\Omega \times [\tau, 3 \tau])$. Iterating this argument we get
	\begin{equation}\label{w L infty L infty}
		w \in L^{\infty}(\Omega \times [\tau, \infty)),\quad\forall \, \tau>0.
	\end{equation}
	
\subsection*{Separation Property.} We are finally in the position to prove the separation property \eqref{separation property} for $\vphi$. We recall that $\vphi$ solves
	\begin{equation}
    \label{NP for vphi}
		-\Delta \vphi + F'(\vphi) = f, \quad \partial_{\n} \vphi|_{\partial \Omega}=0,
	\end{equation}
	with $f = \frac{1}{\ve}( w - \vphi+ \ve \theta_0 \vphi)$. Due to \eqref{w L infty L infty} and since $\|\vphi\|_{L^\infty(\Omega)}\leq 1$, it holds that $f \in L^{\infty}(\Omega \times [\tau, \infty))$, thus, applying \cite[Lemma A.2]{CG}, we deduce that
\begin{equation*}
		\norm{F'(\vphi)}_{L^{\infty}(\Omega \times [\tau, \infty))} \le C.
\end{equation*}
Given that $F'$ diverges at $\pm 1$ and $\vphi \in C(\overline{\Omega} \times [0, \infty))$, we infer the existence of $\delta>0$ that satisfies \eqref{separation property}.
\end{proof}

If the initial datum is more regular, namely $w_0 \in H^1(\Omega)$, we obtain further regularity on the solution.
Therefore, Theorem \ref{ws thm phi and mu} - part (C) can be reformulated in the following way.
\begin{theorem}\label{ss thm}
Under the assumptions of Theorem \ref{ws thm phi and mu} - part (C), let  $w_0= \vphi_0 + \ve \mu_0 \in H^1(\Omega)$ with $\overline{w_0}=m \in \R$. Then, the unique global solution $(w , \vphi)$  to \eqref{w problem}--\eqref{w problem bc i condition} provided by Theorem \ref{ws thm} is strong in the following sense:
\begin{alignat*}{2}
&w \in L^{\infty}(0, \infty; H^1(\Omega)) \cap L^2_{\rm uloc}([0, \infty) ; H^2(\Omega)), \quad \pt w \in L^2_{\rm uloc}([0, \infty) ; L^2(\Omega)),\\
&\vphi \in L^{\infty}(0, \infty; W^{2,p}(\Omega)), \quad \pt \vphi \in L^\infty(0, \infty; H^1(\Omega)),\\
&F'(\vphi)  \in L^{\infty}(0, \infty; L^{p}(\Omega)),
\end{alignat*}
where $p=6$ if $d=3$ and any $2 \le p <\infty$ if $d=2$.
\end{theorem}
\begin{proof}
We will now perform some estimates directly on the weak solution $(w, \vphi)$ of \eqref{w problem}--\eqref{w problem 2}. Those estimates are only formal but can be justified in a suitable approximation scheme. Testing \eqref{w problem} by $- \Delta w$ yields
\begin{equation*}
\frac{\ve}{2} \ddt \norm{\nabla w}_{L^2(\Omega)}^2 
+ \norm{\sqrt{b(\vphi)} \Delta w}_{L^2(\Omega)}^2 
+ \int_{\Omega} b'(\vphi) \nabla \vphi \cdot \nabla w \, \Delta w \,\d x = \int_{\Omega} \div(b(\vphi) \nabla \vphi) \Delta w \,\d x.
\end{equation*}
Applying H{\"o}lder and Young's inequalities together with \eqref{phi L infty H2 unif bound}, we get
\begin{align*}
 \int_{\Omega} b'(\vphi) \nabla \vphi \cdot \nabla w \, \Delta w \,\d x &\le C \norm{\nabla \vphi}_{L^6(\Omega)}\norm{\nabla w}_{L^3(\Omega)} \norm{\Delta w}_{L^2(\Omega)} \notag \\
&\le C \norm{\nabla w}_{L^2(\Omega)}^\frac12 \norm{\Delta w}_{L^2(\Omega)}^\frac32  
\le \frac{b_m}{4} \norm{\Delta w}^2_{L^2(\Omega)} + C \norm{\nabla w}^2_{L^2(\Omega)}
\end{align*}
and
\begin{align*}
\int_{\Omega} \div(b(\vphi) \nabla \vphi) \Delta w \,\d x &\le \frac{b_m}{4} \norm{\Delta w}^2_{L^2(\Omega)} + C\norm{\div(b(\vphi) \nabla \vphi)}^2_{L^2(\Omega)}\nonumber\\
&\le \frac{b_m}{4} \norm{\Delta w}^2_{L^2(\Omega)} 
+ C\left( \norm{b'(\vphi) |\nabla \vphi|^2}_{L^2(\Omega)}^2 + \norm{b(\vphi) \Delta \vphi}_{L^2(\Omega)}^2 \right) \nonumber\\
&\le \frac{b_m}{4} \norm{\Delta w}^2_{L^2(\Omega)} + C \norm{\vphi}_{H^2(\Omega)}^4.
\end{align*}
Thus, exploiting \eqref{phi L infty H2 unif bound}, it holds
\begin{equation}
	\label{diff inq w H1-H2}
	\ve \ddt \norm{\nabla w}_{L^2(\Omega)}^2  + b_m \norm{\Delta w}^2_{L^2(\Omega)} \le C \left( \norm{\nabla w}_{L^2(\Omega)}^2+1\right).
\end{equation}
Therefore, an application of the standard Gr{\"o}nwall's Lemma yields
\begin{equation}
	\label{w L infty H1 unif first part}
	\|\nabla w\|_{L^\infty(0,1;L^2(\Omega))} \le C.
\end{equation}
In order to derive a global in time uniform bound, recalling that $w \in L^2_{\rm uloc}([0, \infty); H^1(\Omega))$, we can apply \cite[Chapter III, Lemma 1.1]{Temam system} with $t_0=0, \ r=1, \ f = \norm{\nabla w}_{L^2(\Omega)}^2$, $g=h=C$ and deduce
\begin{equation}
	\label{w L infty H1 unif control second step}
	\|\nabla w\|_{L^\infty(1,\infty;L^2(\Omega))} \le C.
\end{equation}
Hence, from \eqref{w L infty L2 uniform bound}, \eqref{w L infty H1 unif first part} and \eqref{w L infty H1 unif control second step}  we infer
\begin{equation}
	\label{w L infty H1 unif bound}
	\norm{w}_{L^\infty(0, \infty; H^1(\Omega))} \le C.
\end{equation}
As a consequence, integrating in time \eqref{diff inq w H1-H2} we learn
\begin{align*}
	b_m\int_t^{t+1} \norm{\Delta w(s)}_{L^2(\Omega)}^2 \,\d s &\le \int_t^{t+1}C \left( \norm{\nabla w(s)}_{L^2(\Omega)}^2+1\right) \,\d s + \ve\norm{\nabla w(t)}_{L^2(\Omega)}^2 \nonumber \\
	&\le C \left( \norm{ w}_{L^\infty(0, \infty; H^1(\Omega))}^2+1\right),
\end{align*}
which, together with $w \in L^2_{\rm uloc}([0, \infty); H^1(\Omega))$, implies $w \in L^2_{\rm uloc}([0, \infty); H^2(\Omega))$.
By comparison in \eqref{w problem}, we get
\begin{equation}\label{pt w L2 uloc L2}
	\pt w \in L^2_{\rm uloc}([0, \infty); L^2(\Omega)).
\end{equation}
Besides, in light of \eqref{w problem}, we know that 
$$
\norm{\pt w}_{H^{-1}_{(0)}(\Omega)}\le \frac{C}{\ve}\big(\norm{\nabla w}_{L^2(\Omega)}+\norm{\nabla\vphi}_{L^2(\Omega)}\big).
$$
Thus, combining the above inequality with \eqref{pt vhi H1 abd pt w dual}, we infer from \eqref{phi L infty H1 uniform bound} and \eqref{w L infty H1 unif bound} that
\begin{equation}
	\label{pt phi L2 uloc H1}
	\pt \vphi \in L^\infty(0, \infty; H^1(\Omega)).
\end{equation}
Lastly, recalling that $\vphi$ solves the homogeneous Neumann problem \eqref{NP for vphi}
with $f = \frac{1}{\ve}( w - \vphi+ \ve \theta_0 \vphi) \in L^{\infty}(0, \infty; H^1(\Omega))$, we can apply \cite[Lemma A.4]{CG} to deduce
\begin{equation*}
\norm{\vphi}_{L^\infty(0,\infty; W^{2,p}(\Omega))} + \norm{F'(\vphi)}_{L^\infty(0,\infty; L^p(\Omega))}\le C,
\end{equation*}
where $p=6$ if $d=3$ and any $2 \le p <\infty$ if $d=2$.
\end{proof}

\section{Convergence to the Cahn--Hilliard model}
\label{S-convergenza}

In this section we prove the first part of Theorem \ref{th-singular-limit}. To this aim,
let $\left(\vphi_0, \mu_0\right) \in H^2(\Omega)\times L^2(\Omega)$ satisfy the assumptions of Theorem \ref{ws thm phi and mu} - part (A)
and for any $\ve\in (0,\ve_0)$ let $\vphi^\ve$ be the weak solution to \eqref{CH1-p}--\eqref{nCH-mu} provided by Theorem \ref{ws thm phi and mu} with corresponding chemical potential $\mu^\ve$.

\subsection{Proof of Theorem \ref{th-singular-limit}, Part (A)}
We begin by showing some $\ve$-independent estimates which will be instrumental to study the convergence of the solutions to the $\ve$-perturbed problem
to the solutions of the Cahn-Hilliard system, formally obtained by setting $\ve=0$.
Throughout this section, $C>0$ denotes a generic constant, independent of $\ve$, whose value may change from line to line.

\smallskip
\noindent{\bf Step 1. $\ve$-independent estimates.} 
By the proof of Theorem \ref{ws thm phi and mu}, the following estimates hold true for any
$\ve\in (0,\ve_0)$ (see \eqref{mu L infty L2 uniform bounds}-\eqref{phi L infty H1 uniform bound}):
\begin{align}
\label{uno}
&\norm{\vphi^\ve}_{L^\infty(0,\infty; H^1(\Omega))}+\norm{\sqrt{\ve}\mu^\ve}_{L^\infty(0,\infty; L^2(\Omega))}+
\norm{\nabla \mu^\ve}_{L^2(0, \infty; L^2(\Omega))}\leq C,
\end{align}
We claim that
\begin{equation}
\label{uno-due}
\norm{\mu^\ve}_{L^2_{\rm uloc}([0, \infty); H^1(\Omega))} \le C.
\end{equation}
A key point in this direction is the study of $\overline{\vphi^\ve}$, which is not conserved during the evolution. Nonetheless, we have the conservation law (see \eqref{conservation})
	$$\overline{(\vphi^\ve+\ve\mu^\ve)}(t)=\overline{\vphi_0+\ve\mu_0},\quad \forall \,  t\geq 0,$$
	that is
	\begin{equation}
		\label{emedi}
		\overline{\vphi^\ve} (t)=-\ve \overline{\mu^\ve}(t)
		+\overline{\vphi_0}+\ve\overline{\mu_0},\quad \forall \, t\geq 0.
	\end{equation}
	By the uniform estimate $\sqrt{\ve}\mu^\ve\in L^\infty(0,\infty; L^2(\Omega))$, we deduce that
	$$|\sqrt{\ve}\ \overline{\mu^\ve}(t)| \leq C\int_\Omega \sqrt{\ve}|\mu^\ve(t)|\, \d x \leq C\left(\int_\Omega \ve |\mu^\ve(t)|^2\, \d x \right)^\frac12\leq C.$$
	Then, it is easily seen that
	\begin{equation}
		\label{media phi}
		|\overline{\vphi^\ve}(t)-\overline{\vphi_0}|\leq C\sqrt{\ve},\quad \forall \, t\geq 0,
	\end{equation}
	so $\overline{\vphi^\ve}\to \overline{\vphi_0}$
	uniformly in $[0,\infty)$ as $\ve\to  0$.
	Now observe that, since $|\overline{\vphi_0}|<1$, then $\overline{\vphi_0}\in [-1+2\widetilde{\delta},1-2\widetilde{\delta}]$ for some $\widetilde{\delta}>0$ independent of $\ve$.
	Therefore, there exists $\widetilde\ve>0$ (we assume that $\widetilde\ve<\ve_0$), depending on $\widetilde{\delta}$, such that
	\begin{equation}
		\label{media}
		\overline{\vphi^\ve}(t) \in [-1+\widetilde{\delta},1-\widetilde{\delta} ],\quad \forall \, t\geq 0, \ \forall \, \ve\in (0,\widetilde\ve).
	\end{equation}
In the following, we assume that $\ve\in (0,\widetilde\ve)$. Owing to \eqref{media}, a classical argument first devised in \cite{MZ} (see also \cite[Ch. 4]{MB2019})
allows us to estimate the $L^1$-norm of $F'(\vphi^\ve)$, by means of the following property of the singular potential:
$$
\int_\Omega| F'(u)|\, \d x\leq \widehat{C}\left|\int_\Omega F'(u)(u-\overline{u} )\, \d x\right|+\widehat{C},
$$
where $\widehat{C}$ goes to $+\infty$ as $|\overline{u}|\to 1$.
Accordingly, in light of \eqref{media} we find $\widehat{C}>0$ (depending on $\widetilde{\delta}$, but not on $\ve$) such that, for all $t\geq 0$
	\begin{equation}
		\label{mzz}
		\int_\Omega| F'(\vphi^\ve)|\, \d x\leq \widehat{C}\left|\int_\Omega F'(\vphi^\ve)(\vphi^\ve-\overline{\vphi^\ve})\, \d x\right|+ \widehat{C}.
	\end{equation}
Now, testing the equation for $\mu^\ve$ by $\vphi^\ve-\overline{\vphi^\ve}$ we find
$$(\mu^\ve,\vphi^\ve-\overline{\vphi^\ve})=(-\Delta\vphi^\ve,\vphi^\ve-\overline{\vphi^\ve})+(\Psi'(\vphi^\ve),\vphi^\ve-\overline{\vphi^\ve}),$$
which gives
$$\|\nabla \vphi^\ve\|^2_{L^2(\Omega)}+(F'(\vphi^\ve),\vphi^\ve-\overline{\vphi^\ve})=\theta_0(\vphi^\ve,\vphi^\ve-\overline{\vphi^\ve})+(\mu^\ve-\overline{\mu^\ve},\vphi^\ve-\overline{\vphi^\ve}).$$
As a result, recalling that $\|\vphi^\ve\|_{L^\infty(\Omega)}\leq 1$ and $\nabla\vphi^\ve\in L^\infty(0,\infty; L^2(\Omega))$ uniformly in $\ve$, we get
\begin{align}
\label{mzz0}
    |(F'(\vphi^\ve),\vphi^\ve-\overline{\vphi^\ve})|&\leq \|\nabla \vphi^\ve\|^2_{L^2(\Omega)}+\theta_0|(\vphi^\ve,\vphi^\ve-\overline{\vphi^\ve})|+
|(\mu^\ve-\overline{\mu^\ve},\vphi^\ve-\overline{\vphi^\ve})|\nonumber \\
&\leq C \left( 1+\|\nabla \mu^\ve\|_{L^2(\Omega)} \right).
\end{align}
Combining \eqref{mzz} and \eqref{mzz0}, we thus obtain
$$
\int_\Omega| F'(\varphi^\ve)|\, \d x 	\leq 
C \left( 1+\|\nabla \mu^\ve\|_{L^2(\Omega)} \right).
$$
Since $\overline{\mu^\ve} =|\Omega|^{-1}\int_\Omega F'(\vphi^\ve)\,\d x -\theta_0\overline{\vphi^\ve}$, owing to \eqref{media}, we arrive at
\begin{equation}
\label{mediamu}
|\overline{\mu^\ve}|\leq C\left(1+ \int_\Omega| F'(\varphi^\ve)|\,\d x\right)\leq C \left( 1+\|\nabla \mu^\ve\|_{L^2(\Omega)} \right)\!,
\end{equation}
which, together with the uniform bound $\norm{\nabla \mu^\ve}_{L^2(0, \infty; L^2(\Omega))} \le C$, concludes the proof of \eqref{uno-due}. Then, as a consequence of \eqref{uno} and \eqref{uno-due}, we also find the uniform control
\begin{equation}
\label{tre}
\norm{\vphi^\ve + \ve \mu^\ve}_{L^2_{\rm uloc}([0, \infty); H^1(\Omega))}\leq C.
\end{equation}

Next, observing that $\vphi^\ve$ solves the singular elliptic problem
\begin{equation}
\label{NEUMANNSING}
\begin{cases}
-\Delta \vphi^\ve+F'(\vphi^\ve)=f,\quad &\text{ in }\Omega,\\
\partial_\n \vphi^\ve=0, \quad &\text{ on }\partial \Omega,
\end{cases}
\end{equation}
with $f=\mu^\ve+\theta_0\vphi^\ve$,
an application of \cite[Lemma A.3]{CG} and \eqref{uno} yields
\begin{align}
\label{phi h2 nabla mu}
\norm{\Delta \vphi^\ve}_{L^2(\Omega)}
\le \norm{\nabla \mu^\ve+ \theta_0 \nabla \vphi^\ve}_{L^2(\Omega)}^\frac12 \norm{\nabla \vphi^\ve}_{L^2(\Omega)}^\frac12
\le C \left( 1+\|\nabla \mu^\ve\|_{L^2(\Omega)}^\frac12 \right)\!,
\end{align}
which, together with \eqref{uno-due}, implies
	\begin{align}
		\label{convergence vphi ve H2}
		\norm{\vphi^\ve}_{L^4_{\rm uloc}([0, \infty); H^2(\Omega))} \le C.
	\end{align}
Lastly, from \cite[Lemma A.1]{CG} we infer that
	\begin{equation*}
		\norm{F'(\vphi^\ve)}_{L^2(\Omega)} \le \norm{\mu^\ve+ \theta_0\vphi^\ve}_{L^2(\Omega)} \le C \left( 1+\|\nabla \mu^\ve\|_{L^2(\Omega)} \right),
	\end{equation*}
	which, by \eqref{uno-due} entails that
	\begin{equation}
\label{quattro}
		\norm{F'(\vphi^\ve)}_{L^2_{\rm uloc}([0, \infty); L^2(\Omega))} \le C.
	\end{equation}

\subsection*{Step 2. Passage to the limit $\ve\to 0$.}
 Let $T>0$ be fixed. By the $\ve$-uniform estimates above, together with \eqref{pt w L2 dual uniform bound},
 we infer that, up to subsequences, the following convergences hold true as $\ve \to 0$:
\begin{alignat}{2}
\vphi^\ve &\rightharpoonup \vphi^\ast \quad &&\text{weakly in}\  L^4(0,T;H^2(\Omega)), 
\label{convergence phi with perturbed energy equality}\\
\vphi^\ve+ \ve \mu^\ve &\rightharpoonup \widetilde{w} \quad &&\text{weakly in}\ L^2(0, T; H^1(\Omega)), 
\\
\pt (\vphi^\ve+ \ve \mu^\ve) &\rightharpoonup z \quad &&\text{weakly in}\ L^2(0,T; (H^1(\Omega))'), 
\label{weak convergence pt w}\\
b(\vphi^\ve)\nabla \mu^\ve &\rightharpoonup \omega \quad &&\text{weakly in}\ L^2(\Omega \times(0, T)),
\end{alignat}
for some functions $\widetilde{w},z, \omega$ to be identified.
Besides, since $\mu^\ve \in C([0,T]; L^2(\Omega))$ for any $T>0$ and any $\ve >0$, we have 
\begin{equation*}
	\sup_{t \in [0,T]}\norm{\ve \mu^\ve(t)}_{L^2(\Omega)} = \sqrt{\ve} \sup_{t \in [0,T]} \norm{\sqrt{\ve} \mu^\ve(t)}_{L^2(\Omega)} \le C \sqrt{\ve},
\end{equation*}
which yields
\begin{equation}
	\label{convergence mu with perturbed equlity}
	\ve \mu^\ve \to 0 \quad \text{strongly in}\ C([0, T]; L^2(\Omega)), \quad  \text{as}\ \ve \to 0.
\end{equation}
Since $\vphi^\ve+\ve\mu^\ve$ is bounded in $L^\infty(0,T;L^2(\Omega))\cap L^2(0,T;H^1(\Omega))$ and its time derivative in $L^2(0,T;(H^1(\Omega))')$, an application of the Aubin--Lions lemma (see \cite[Chapter III, Theorem 2.1]{Temam2001} and \cite[Corollary 4]{Simon1987}) yields
\begin{equation}
	\label{w converges contin}
	\vphi^\ve + \ve \mu^\ve \to \widetilde{w} \quad \text{strongly in} \ L^2(0, T; L^2(\Omega)) \cap C([0,T];(H^1(\Omega))'), \quad \ \text{as}\ \ve \to 0.
\end{equation}
This allows us to identify $\widetilde{w}$. Indeed, we observe that
\begin{equation}
	\label{phi strong convergence}
	\vphi^\ve  = (\vphi^\ve + \ve \mu^\ve )- \ve \mu^\ve  \to \widetilde{w} \quad \text{strongly in} \ L^2(0, T; L^2(\Omega)), \quad  \text{as}\ \ve \to 0.
\end{equation}
Hence, owing to \eqref{convergence phi with perturbed energy equality} and \eqref{phi strong convergence}, we deduce that
\begin{equation}
	\label{limit w ve}
	\widetilde{w} = \vphi^\ast
\end{equation}
and
\begin{equation}
	\label{a.e. phi convergence}
	\vphi^\ve\to \vphi^\ast \quad \text{a.e. in }\Omega\times (0, T), \quad \text{as}\ \ve \to 0.
\end{equation}
By \eqref{convergence mu with perturbed equlity} and \eqref{w converges contin}, we also infer that $\vphi^\ve\to\vphi^0$ in $C([0,T];(H^1(\Omega))')$. Then, by \eqref{Hilbert interpolation inequality} and \eqref{uno},
$$
\sup_{t\in[0,T]}\norm{\vphi^\ve(t)-\vphi^\ast(t)}_{L^2(\Omega)}\le C\sup_{t\in[0,T]}\norm{\vphi^\ve(t)-\vphi^\ast(t)}_{(H^1(\Omega))'}^{\frac12}
\to 0,
$$
that is, $\vphi^\ve\to\vphi^\ast$ in $C([0,T];L^2(\Omega))$. Finally, the interpolation inequality 
$\norm{v}_{H^{2-\alpha}(\Omega)}\le C\norm{v}_{H^2(\Omega)}^{1-\frac\alpha2}\norm{v}_{L^2(\Omega)}^{\frac\alpha2}$, for $\alpha\in(0,2]$, 
and \eqref{convergence vphi ve H2} give
$$
\int_0^T\norm{\vphi^\ve(t)-\vphi^\ast}(t)_{H^{2-\alpha}(\Omega)}^4\,\d t\le C\sup_{t\in[0,T]}\norm{\vphi^\ve(t)-\vphi^\ast(t)}_{L^2(\Omega)}^{2\alpha}\int_0^T\norm{\vphi^\ve(t)-\vphi^\ast(t)}_{H^2(\Omega)}^{4-2\alpha}\,\d t
\to 0,
$$
which proves the first convergence in Theorem \ref{th-singular-limit} - (A).

Next, we identify $z$. Given \eqref{weak convergence pt w} and \eqref{limit w ve}, for any $\eta \in C^\infty_{c}(0, T)$ and for any $v \in H^1(\Omega)$, for $ \ve \to 0$, it follows that
\begin{align*}
	&\int_0^T \eta(s)\l \pt (\vphi^\ve+\ve\mu^\ve)(s), v \r \,\d s = - \int_0^T \pt \eta(s) ((\vphi^\ve+\ve\mu^\ve)(s), v)\, \d s \to - \int_0^T \pt \eta(s) (\vphi^\ast (s), v) \,\d s, 
\end{align*}
and, at the same time,
\begin{align*}
	&\int_0^T \eta(s)\l \pt (\vphi^\ve+\ve\mu^\ve)(s), v \r \,\d s \to \int_0^T \eta(s) \l z(s), v \r \,\d s,
\end{align*}
which imply that
\begin{equation}
	\label{phi limit has a time derivative}
	\pt (\vphi^\ve+\ve\mu^\ve) \rightharpoonup \pt \vphi^\ast \quad \text{weakly in}\ L^2(0, T; (H^1(\Omega))'),\quad \text{as}\ \ve \to 0.
\end{equation}
\smallskip
We now use a classical argument to prove $|\vphi^\ast| < 1$ a.e. in $\Omega\times (0, T)$. For all $n \in \mathbb{N}$, we introduce the sets
\begin{equation*}
	\begin{split}
		E_n = \left\{(x,t) \in \Omega \times (0,T): |\vphi^\ast(x,t) |> 1- \frac{1}{n} \right\}, \\
		E_n^{\ve} = \left\{(x,t) \in \Omega \times (0,T): |\vphi^\ve (x,t)| > 1- \frac{1}{n} \right\}.
	\end{split}
\end{equation*}
For any fixed $n$, an application of Fatou's Lemma together with the pointwise convergence of $\vphi^\ve $ yields
\begin{equation*}
	| E_n | \le \liminf_{\ve \to 0^+} |E_n^{\ve}| .
\end{equation*}
Observe that, in light of \eqref{quattro}, we have $\norm{F'(\vphi^\ve)}_{L^1(\Omega\times (0,T))} \le C$. Besides,
\begin{equation*}
	\begin{split}
		\norm{F' \left(\vphi^\ve\right)}_{L^1(\Omega \times (0, T))}& \ge \int_{E_n^{\ve}} |F'\left(\vphi^\ve\right)|\, \d x \d t \\
		&\ge \min\left\{F'\left(1-\frac{1}{n}\right), - F'\left(-1+\frac{1}{n}\right)\right\} |E_n^\ve|.
	\end{split}
\end{equation*}
Therefore, we obtain
\begin{equation*}
|E_n| \le \frac{C}{\min\{F'\left(1-\frac{1}{n}\right), - F'\left(-1+\frac{1}{n}\right)\}}.
\end{equation*}
Taking the limit for $n \to \infty$, we deduce that
\begin{equation*}
|\lbrace (x,t) \in \Omega \times (0,T): |\vphi^\ast(x,t)| \ge 1 \rbrace |=0,
\end{equation*}
which gives the desired conclusion.

Let us now observe that, since $F' \in C^0(-1,1)$, by the pointwise convergence \eqref{a.e. phi convergence}, it holds
$F'(\vphi^\ve) \to F'(\vphi^\ast)$ a.e in $\Omega\times (0,T)$. Then, thanks to \eqref{quattro}, the weak dominated convergence Theorem (cfr. \cite[Lemma 1.3]{lions1969}) implies
\begin{equation}
	\label{limit F'}
	F'(\vphi^\ve) \rightharpoonup F'(\vphi^\ast) \quad \text{weakly in}\ L^2(\Omega\times (0,T)).
\end{equation}
Besides, since $\Delta \vphi^\ve\rightharpoonup \Delta \vphi^\ast$ in $L^4(0,T;L^2(\Omega))$ by \eqref{convergence phi with perturbed energy equality}, we have
\begin{equation*}
	\mu^\ve= -\Delta \vphi^\ve + F'(\vphi^\ve)- \theta_0\vphi^\ve \rightharpoonup -\Delta \vphi^\ast + F'(\vphi^\ast) - \theta_0 \vphi^\ast \quad \text{weakly in}\ L^2(\Omega\times(0,T)),\quad \text{as} \ \ve \to 0.
\end{equation*}
We set
\begin{equation}
	\label{mu 0 non deg }
	\mu^\ast =  -\Delta \vphi^\ast + F'(\vphi^\ast) - \theta_0\vphi^\ast \quad \text{a.e. in}\ \Omega\times(0,T) .
\end{equation}
Note that, in light of the uniform bound  of $\mu^\ve$ in $L^2_{\rm uloc}([0,\infty); H^1(\Omega))$, we also have
\begin{equation}
\label{mu weak convergence}
\nabla \mu^\ve \rightharpoonup \nabla \mu^\ast \quad \text{weakly in}\ L^2(\Omega\times(0,T) ).
\end{equation}
Lastly, we show that $\omega = b(\vphi^\ast)\nabla \mu^\ast$. 
To this aim, we observe that \eqref{a.e. phi convergence} and the continuity of $b$ on $[-1,1]$ ensure that $b(\vphi^\ve) \to b(\vphi^\ast)$ a.e. in $\Omega\times (0, T)$. Being $b$ bounded, an application of Lebesgue's Theorem yields
\begin{equation}
	\label{mobility convergence}
	b(\vphi^\ve) \to b(\vphi^\ast) \quad \text{strongly in}\ L^p(\Omega\times (0, T)), \ \forall \, p< \infty, \quad \text{as}\ \ve \to 0.
\end{equation}
Hence, we deduce by \eqref{mu weak convergence}
that
\begin{equation}
	\label{convergence rhs}
	b(\vphi^\ve) \nabla \mu^\ve \rightharpoonup b(\vphi^\ast)\nabla \mu^\ast \quad \text{weakly in}\ L^2(\Omega\times(0,T) ), \quad \text{as}\ \ve \to 0.
\end{equation}

Next, we consider $\eta \in C^1_{c}(\overline\Omega \times [0, \infty))$ as a test function in \eqref{weak formulation perturbed Cahn--Hilliard equation}, which gives
$$ 
\int_\Omega \l \pt \vphi^\ve + \ve \pt \mu^\ve \r \eta \,\d x  = - \int_\Omega b(\vphi^\ve) \nabla \mu^\ve \cdot \nabla \eta \, \d x, 
$$
for almost every $t \in (0,\infty)$. An integration in time with $T \in (0,\infty)$ gives
\begin{align}
	\label{very weak perturbed Cahn Hilliard}
	\int_0^T &\int_\Omega \vphi^\ve \pt \eta  \,\d x \, \d t + \int_0^T \int_\Omega \ve \mu^\ve \pt \eta \, \d x \, \d t \nonumber \\
	&= \int_0^T \int_\Omega b(\vphi^\ve)\nabla \mu^\ve \cdot \nabla \eta  \,\d x \, \d t + \int_\Omega \left(\vphi^\ve(T)+ \ve \mu^\ve(T) \right) \eta(T)\, \d x - \int_{\Omega} \left(\vphi_0 + \ve \mu_0 \right) \eta(0)\, \d x.
\end{align}
Exploiting \eqref{convergence phi with perturbed energy equality}, \eqref{convergence mu with perturbed equlity}, \eqref{w converges contin}, \eqref{limit w ve}, \eqref{convergence rhs},  we can pass to the limit in \eqref{very weak perturbed Cahn Hilliard} (for the term at time $T$ we use the convergence in $C([0,T];(H^1(\Omega))')$ in \eqref{w converges contin}) and infer that $\left(\vphi^\ast, \mu^\ast\right)$ solves
\begin{equation}
    \label{very weak Cahn Hilliard}
    \int_0^T \int_\Omega \vphi^\ast \pt \eta  \,\d x \, \d t = \int_0^T \int_\Omega b(\vphi^\ast)\nabla \mu^\ast \cdot \nabla \eta  \,\d x \, \d t + \int_\Omega \vphi^\ast(T) \eta(T) \,\d x
	- \int_{\Omega} \vphi^\ast \eta(0) \,\d x,
\end{equation}
for all $ \eta \in C^1_{c} (\overline\Omega \times [0, \infty))$, $T\in (0,\infty)$.
Hence, recalling \eqref{phi limit has a time derivative}, \eqref{mu 0 non deg } and exploiting \cite[Lemma 7.4]{Robinson system}, we deduce that $\left(\vphi^\ast, \mu^\ast \right)$ solves \eqref{CH1}--\eqref{CH2} for any $T>0$ in the sense of Definition \ref{def-CH-weak}. 
\smallskip


\section{Rate of convergence in two dimensions}
\label{Sec-2D}

In this section, we restrict our analysis to the two-dimensional case and focus on proving the second part of Theorem \ref{th-singular-limit}. Our first goal is to derive additional $\ve$-independent estimates in stronger norms, which will play a key role in establishing the rate of convergence of the solutions as $\ve\to 0$.
\begin{theorem}
	\label{t-uniform}
	Let $d=2$. Let $\left(\vphi_0,\mu_0\right)\in H^2(\Omega)\times H^1(\Omega)$ and assume the hypothesis of Theorem \ref{th-singular-limit} - (B). Then, there exists $\ve_1 \in (0, \ve_0)$ such that,
	for every $\ve\in (0,\ve_1)$, the corresponding strong solution $\vphi^\ve$ to \eqref{CH1-p}--\eqref{CH2-p} subject to \eqref{nCH-mu}, with corresponding $\mu^\ve$, has the following properties:
 \begin{align*}
		\norm{\vphi^\ve}_{L^\infty(0,\infty; H^2(\Omega))} &+  \norm{\pt \vphi^\ve}_{L^2(0,\infty; H^1(\Omega))} \le C,\\
		\norm{\mu^\ve}_{L^\infty(0,\infty; H^1(\Omega))} &
+ \norm{\sqrt{\ve}\pt \mu^\ve}_{L^2(0,\infty; L^2(\Omega))}+\norm{F'(\varphi^\ve)}_{L^\infty(0,\infty; L^2(\Omega))} \le C,
	\end{align*}
	where $C>0$ is independent of $\ve$.
\end{theorem}
%
\begin{proof} We start by exhibiting uniform bounds for the norms of $\mu^\ve$. To this aim, let us first observe that by \eqref{mediamu}
\begin{equation}
\label{H2mu}
\|\mu^\ve\|_{H^2(\Omega)}\leq C(1+ \norm{\nabla\mu^\ve}_{L^2(\Omega)}+\norm{\Delta\mu^\ve}_{L^2(\Omega)}).
\end{equation}
Testing \eqref{CH1-p} by $\Delta \mu^\ve$, we have
\begin{equation}
\label{uffa}
b_m\|\Delta\mu^\ve\|_{L^2(\Omega)}^2\leq \int_\Omega b(\vphi^\ve)|\Delta\mu^\ve|^2 \,\d x=(\pt \vphi^\ve,\Delta \mu^\ve)+\ve(\pt \mu^\ve,\Delta \mu^\ve)
-\int_\Omega b'(\vphi^\ve)\nabla\vphi^\ve \cdot \nabla \mu^\ve\Delta \mu^\ve \,\d x.
\end{equation}
Exploiting Ladyzhenskaya's inequality $\|\nabla u\|_{L^4(\Omega)}^2\leq C \|\nabla u\|_{L^2(\Omega)}\|u\|_{H^2(\Omega)}$, the uniform bounds \eqref{uno} and \eqref{phi h2 nabla mu}, we obtain
\begin{align*}
	\int_\Omega &b'(\vphi^\ve)\nabla\vphi^\ve \cdot \nabla \mu^\ve\Delta \mu^\ve \,\d x \\
    &\leq C||\nabla\vphi^\ve\|_{L^4(\Omega)}\|\nabla\mu^\ve\|_{L^4(\Omega)}\|\Delta \mu^\ve\|_{L^2(\Omega)}\\
	&\leq C\|\nabla \vphi^\ve\|_{L^2(\Omega)}^\frac12\|\vphi^\ve\|_{H^2(\Omega)}^\frac12\|\nabla\mu^\ve\|_{L^2(\Omega)}^\frac12\|\mu^\ve\|_{H^2(\Omega)}^\frac12\norm{\Delta\mu^\ve}_{L^2(\Omega)}\\
	&\le C\left(1+\norm{\Delta \vphi^\ve}_{L^2(\Omega)}^\frac12\right)\norm{\nabla \mu^\ve}_{L^2(\Omega)}^\frac12
\left(1+ \norm{\nabla\mu^\ve}_{L^2(\Omega)}^\frac12+\norm{\Delta\mu^\ve}_{L^2(\Omega)}^\frac12\right)\norm{\Delta\mu^\ve}_{L^2(\Omega)}\\
	&\leq C \left(\|\nabla \mu^\ve\|_{L^2(\Omega)}^3 + \norm{\nabla\mu^\ve}_{L^2(\Omega)}^\frac32+ \norm{\nabla\mu^\ve}_{L^2(\Omega)}^2+ \norm{\nabla\mu^\ve}_{L^2(\Omega)}^\frac52+\norm{\nabla\mu^\ve}_{L^2(\Omega)}\right)+\dfrac{b_m}2\|\Delta \mu^\ve\|_{L^2(\Omega)}^2.
\end{align*}
Thus, collecting \eqref{H2mu} and \eqref{uffa}, we find
\begin{align}
	\label{deltamu}
\|\mu^\ve\|_{H^2(\Omega)}^2&\leq C \|\pt \vphi^\ve\|_{L^2(\Omega)}^2+C\ve^2\|\pt \mu^\ve\|_{L^2(\Omega)}^2
	+C\left(1+\|\nabla \mu^\ve\|_{L^2(\Omega)}^3 \right).
\end{align}
At this point, we look for uniform estimates of $\nabla \mu^\ve$, $\pt\vphi^\ve$ and $\pt\mu^\ve$.
Testing \eqref{CH1-p} by $\pt\mu^\ve$ and recalling \eqref{integration by part formula}, we get
$$-(\pt\vphi^\ve, \pt\mu^\ve)+\frac12\int_\Omega b'(\vphi^\ve)\pt\vphi^\ve|\nabla \mu^\ve|^2\,\d x=\frac12\ddt\int_\Omega b(\vphi^\ve)|\nabla \mu^\ve|^2\,\d x+\ve\|\pt\mu^\ve\|_{L^2(\Omega)}^2.$$
Differentiating \eqref{CH2-p} with respect to time, 
$$
\pt\mu^\ve=-\Delta \pt\vphi^\ve+\Psi''(\vphi^\ve)\pt\vphi^\ve,
$$
and testing the latter by $\pt\vphi^\ve$, we find
$$(\pt\mu^\ve,\pt\vphi^\ve)=\|\nabla \pt \vphi^\ve\|_{L^2(\Omega)}^2+ (\Psi''(\vphi^\ve)\pt \vphi^\ve,\pt \vphi^\ve) .$$
This yields
\begin{equation}
    \label{test con ptmu}
    \frac12\ddt\int_\Omega b(\vphi^\ve)|\nabla \mu^\ve|^2\,\d x+\ve\|\pt \mu^\ve\|^2_{L^2(\Omega)}+\|\nabla \pt \vphi^\ve\|_{L^2(\Omega)}^2\leq \theta_0\|\pt \vphi^\ve\|_{L^2(\Omega)}^2+\frac12\int_\Omega b'(\vphi^\ve)\pt \vphi^\ve|\nabla \mu^\ve|^2\,\d x.
\end{equation}
Adding $\|\pt \vphi^\ve\|_{L^2(\Omega)}^2$ to both sides of \eqref{test con ptmu}, we observe that by the interpolation inequality \eqref{Hilbert interpolation inequality}
$$(1+\theta_0)\|\pt \vphi^\ve\|_{L^2(\Omega)}^2\leq (1+\theta_0)\|\pt \vphi^\ve\|_{(H^1(\Omega))'}\|\pt \vphi^\ve\|_{H^1(\Omega)}\leq \frac12 \|\pt \vphi^\ve\|_{H^1(\Omega)}^2+\frac{(1+\theta_0)^2}{2}\|\pt \vphi^\ve\|_{(H^1(\Omega))'}^2.$$
Besides, since
$\|\pt \vphi^\ve+\ve \pt \mu^\ve\|_{(H^1(\Omega))'}\le b_M\|\nabla \mu^\ve\|_{L^2(\Omega)}$ by \eqref{CH1-p},
we have
$$\|\pt \vphi^\ve\|_{(H^1(\Omega))'}\leq \ve \|\pt \mu^\ve\|_{(H^1(\Omega))'}+C\|\nabla \mu^\ve\|_{L^2(\Omega)}\leq \ve C_1\|\pt \mu^\ve\|_{L^2(\Omega)}+C\|\nabla \mu^\ve\|_{L^2(\Omega)},$$
for some $C_1>0$ given by the embedding $L^2(\Omega)\subset (H^1(\Omega))'$.
Therefore, restricting $\ve$ so that $\ve\in (0,(1+\theta_0)^{-2}C_1^{-2}/2)$, we obtain
\begin{equation}
\ddt\int_\Omega b(\vphi^\ve)|\nabla \mu^\ve|^2\,\d x
+\ve\|\pt \mu^\ve\|_{L^2(\Omega)}^2
+\|\pt \vphi^\ve\|_{H^1(\Omega)}^2
\leq C\|\nabla \mu^\ve\|_{L^2(\Omega)}^2
+\int_\Omega b'(\vphi^\ve)\pt \vphi^\ve|\nabla \mu^\ve|^2\,\d x.
\label{uniffoo}
\end{equation}
We proceed by estimating the last term. By Lemma \ref{L2 product lemma}, 
\begin{align}
\label{P}
	P=\int_\Omega b'(\vphi^\ve)\pt \vphi^\ve|\nabla \mu^\ve|^2\,\d x
	&\leq C\|\nabla\mu^\ve\|_{L^2(\Omega)}\|\pt \vphi^\ve\nabla \mu^\ve\|_{L^2(\Omega)} \notag \\
	&\leq C\|\nabla\mu^\ve\|_{L^2(\Omega)}^2\|\pt \vphi^\ve\|_{H^1(\Omega)}\ln^\frac12\left(C\frac{\|\mu^\ve\|_{H^2(\Omega)}}{\|\nabla \mu^\ve\|_{L^2(\Omega)}}\right) \notag \\
	&\leq \frac14\|\pt \vphi^\ve\|_{H^1(\Omega)}^2+C\|\nabla\mu^\ve\|_{L^2(\Omega)}^4
\ln\left(C\frac{\|\mu^\ve\|_{H^2(\Omega)}}{\|\nabla \mu^\ve\|_{L^2(\Omega)}}\right).
\end{align}
Using the elementary inequalities
\begin{align*}
    x^2\ln\left(C\frac{y}{x}\right)\leq x^2\ln(eCy)+1 \quad \text{and} \quad 
    xy\leq x\ln x+e^y,\quad \forall x>0, y>0,
\end{align*}
we have
\begin{align*}
\|\nabla\mu^\ve\|_{L^2(\Omega)}^2\ln\left(C\frac{\|\mu^\ve\|_{H^2(\Omega)}}{\|\nabla \mu^\ve\|_{L^2(\Omega)}}\right)&\leq \|\nabla\mu^\ve\|_{L^2(\Omega)}^2\ln(eC\|\mu^\ve\|_{H^2(\Omega)})+1\\
	&\leq \|\nabla\mu^\ve\|_{L^2(\Omega)}^2\ln(\|\nabla \mu^\ve\|_{L^2(\Omega)}^2)+eC\|\mu^\ve\|_{H^2(\Omega)}+1.
\end{align*}
Hence,
\begin{align*}
	P\leq \frac14\|\pt \vphi^\ve\|_{H^1(\Omega)}^2+C\|\nabla\mu^\ve\|_{L^2(\Omega)}^2+ C\|\nabla\mu^\ve\|_{L^2(\Omega)}^4\ln(\|\nabla \mu^\ve\|_{L^2(\Omega)}^2)+C\|\mu^\ve\|_{H^2(\Omega)} \|\nabla\mu^\ve\|_{L^2(\Omega)}^2.
\end{align*}
Recalling \eqref{deltamu}, we deduce that
\begin{align*}
	P &\le \frac14\|\pt \vphi^\ve\|_{H^1(\Omega)}^2+C\|\nabla\mu^\ve\|_{L^2(\Omega)}^2+C\|\nabla\mu^\ve\|_{L^2(\Omega)}^4\ln(\|\nabla \mu^\ve\|_{L^2(\Omega)}^2)\nonumber\\
	&\quad +C\left(1+\norm{\pt \vphi^\ve}_{L^2(\Omega)}+ \ve \norm{\pt \mu^\ve}_{L^2(\Omega)} + \norm{\nabla \mu^\ve}^2_{L^2(\Omega)}\right)\norm{\nabla \mu^\ve}^2_{L^2(\Omega)} \nonumber \\
	&\le \frac12 \norm{\pt \vphi^\ve}_{H^1(\Omega)}^2 + C_2 \ve^2 \norm{\pt \mu^\ve}^2_{L^2(\Omega)} +C \norm{\nabla \mu^\ve}^2_{L^2(\Omega)}+ C\|\nabla\mu^\ve\|_{L^2(\Omega)}^4\ln(e+ \|\nabla \mu^\ve\|_{L^2(\Omega)}^2).
\end{align*}
As a result, assuming $\ve<\varepsilon_1$ where
$
\ve_1:=\min\{\widetilde\ve,\,(1+\theta_0)^{-2}C_1^{-2}/2,\,(2C_2)^{-1}\},
$
we  end up with
$$
\ddt\int_\Omega b(\vphi^\ve)|\nabla \mu^\ve|^2\,\d x +\frac\ve2\|\pt \mu^\ve\|_{L^2(\Omega)}^2+\frac12\|\pt \vphi^\ve\|_{H^1(\Omega)}^2\\
\leq C\|\nabla \mu^\ve\|_{L^2(\Omega)}^2+ C\|\nabla\mu^\ve\|_{L^2(\Omega)}^4\ln(e+\|\nabla\mu^\ve\|_{L^2(\Omega)}^2).
$$
Setting
$$
B(t)=\int_\Omega b(\vphi^\ve(t))|\nabla \mu^\ve(t)|^2 \, \d x,
$$
and observing that
$$b_m \norm{\nabla \mu^\ve (t)}^2_{L^2(\Omega)} \le B(t) \le b_M \norm{\nabla \mu^\ve (t)}^2_{L^2(\Omega)},$$
we finally have
\begin{equation}
	\label{uniff}
	\ddt B(t)+ \frac\ve2\|\pt \mu^\ve\|_{L^2(\Omega)}^2+\frac12\|\pt \vphi^\ve\|_{H^1(\Omega)}^2\leq CB(t)+ CB(t)^2\ln(e+B(t)), \quad \text{for a.e. }t\geq 0.
\end{equation}
where $B\in L^1(0,\infty)$ in light of \eqref{uno} and $B(0)\leq b_M\|\nabla\mu_0\|_{L^2(\Omega)}^2 <\infty$ by assumption.
Therefore, by the logarithmic Gr{\"o}nwall Lemma \ref{logaritmic Gronwall}, we conclude that
\begin{equation}
	\label{mutfit1}
	\|\nabla\mu^\ve\|_{L^\infty(0,\infty; L^2(\Omega))}\leq C.
\end{equation}
Finally, an integration of \eqref{uniff} on $(0,+\infty)$ yields
\begin{equation}
	\label{mutfit2}
	\| \pt \vphi^\ve\|_{L^2(0,\infty; H^1(\Omega))}+\|\sqrt{\ve}\pt \mu^\ve\|_{L^2(0,\infty; L^2(\Omega))}\leq C.
\end{equation}
In light of \eqref{mediamu}, we infer from \eqref{mutfit1} that
\begin{equation}
	\label{duez}
	\|\mu^\ve\|_{L^\infty(0,\infty; H^1(\Omega))}\leq C.
\end{equation}
\smallskip
\noindent
Finally, we observe once again that $\vphi^\ve$ solves the elliptic problem \eqref{NEUMANNSING}
with $f= \mu^\ve + \theta_0 \vphi^\ve$.
Owing to \eqref{uno} and \eqref{duez} we learn that
$\|f\|_{L^\infty(0,\infty; H^1(\Omega))}\leq C$, and thus an application of \cite[Lemma A.4]{CG} yields
$$\|\vphi^\ve\|_{L^\infty(0,\infty;W^{2,p}(\Omega))}+\|F'(\vphi^\ve)\|_{L^\infty(0,\infty;L^p(\Omega))}\leq C,$$
for any finite $p\geq 2$ since $d=2$.
The proof is completed.
\end{proof}

\begin{remark}
With the additional bounds obtained in Theorem \ref{t-uniform}, one can pass to the limit $\ve \to 0$ with stronger convergences. In particular
\begin{align*}
\label{convercence with respect to eps 1}
&\vphi^\ve \stackrel{\ast}{\rightharpoonup} \vphi^\ast \quad \text{weak-star in} \ L^{\infty}(0, T; H^2(\Omega)),\qquad \pt \vphi^\ve \rightharpoonup \pt \vphi^\ast \quad\text{weakly in}\ L^2(0, T; H^1(\Omega)),\nonumber \\
	&\mu^\ve \stackrel{\ast}{\rightharpoonup} \mu^\ast \quad \text{weak-star in}\ L^{\infty}(0, T; H^1(\Omega)),\qquad 
	\ve \pt \mu^\ve \to 0 \quad\text{strongly in}\  L^2(0, T; L^2(\Omega)).
\end{align*}
Moreover, exploiting the Aubin-Lions lemma, we infer that $\vphi^\ve \to \vphi^\ast$ in $L^2(0, T; H^1(\Omega)) \cap C([0, T]; L^2(\Omega))$, and therefore $\vphi^\ve(x, t) \to \vphi^\ast(x, t)$ almost everywhere in $\Omega \times (0, T)$. Then, passing to the limit in the weak formulation of \eqref{CH1-p}-\eqref{CH2-p}, we can show that the limit $(\vphi^\ast, \mu^\ast)$ is indeed the (strong) solution to the Cahn--Hilliard system with non-degenerate mobility (cf. \cite{CGGG2025}).
\end{remark}
\begin{remark} Another interesting consequence of the above uniform bounds is that the strict separation property of $\vphi^\ve$ from the pure phases $\pm 1$ does not depend on $\ve$, provided the nonlinearity $F$ complies with the additional assumption 
\begin{equation}
\label{Fsecondo}
F''(s) \le C e^{C|F'(s)|+C}, \quad \forall s \in (-1,1),
\end{equation}
which is in particular satisfied by the logarithmic potential.
Indeed, $\vphi^\ve$ solves the elliptic problem \eqref{NEUMANNSING}
with $f= \mu^\ve + \theta_0 \vphi^\ve$ and $\norm{f}_{L^\infty(0, \infty; H^1(\Omega))}\le C$.
Therefore, owing to \eqref{Fsecondo} and reasoning as in \cite{HW} (see also \cite[Theorem 6.2]{CGM2024}), we deduce that
$$
|\vphi^\ve(x, t)| \le 1- \delta \quad \forall\, (x, t) \in \overline{\Omega} \times [0, \infty),
$$
where $\delta$ is independent of $\varepsilon$.
\end{remark}

\subsection{Proof of Theorem \ref{th-singular-limit} part (B)}
Firstly, we reformulate the $\ve$-problem as
$$\mu^\ve- \overline{\mu^\ve}=-\G_{\varphi^\ve} (\pt\vphi^\ve-\overline{\pt\vphi^\ve})-\ve\G_{\varphi^\ve} (\pt \mu^\ve-\overline{\pt \mu^\ve})$$
while, owing to the mass conservation of $\vphi^\ast$, the Cahn--Hilliard equation can be written as (cf. Definition \ref{def-CH-weak})
\begin{equation}
\label{CHN}
	\mu^\ast - \overline{\mu^\ast }=-\G_{\vphi^\ast } (\pt \vphi^\ast).
\end{equation}
\noindent
Setting $$
M=\mu^\ve-\mu^\ast  \quad\text{and}\quad \Phi=\vphi^\ve-\vphi^\ast,
$$
we observe that
$$\overline{\Phi(t)}=-\ve\overline{\mu^\ve(t)}+\ve\overline{\mu_0},$$
due to the mass conservation of $\vphi^\ve+\ve\mu^\ve$ and $\vphi^0$.
Since $\norm{\mu^\ve}_{L^\infty(0,\infty; H^1(\Omega))} \le C$ by Theorem \ref{t-uniform}, we infer that
\begin{equation}
	\label{mee}
	|\overline{\Phi(t)}|\leq C{\ve}, \quad \forall \, t\geq 0.
\end{equation}
Besides, $(\Phi, M)$ solves the system
\begin{subequations}
	\begin{alignat}{1}
		\label{sistem difference 1}
		&M - \overline{M}=-\G_{\varphi^\ve} (\pt \vphi^\ve-\overline{\pt \vphi^\ve})+\G_{\vphi^\ast} (\pt \vphi^\ast )-\ve\G_{\varphi^\ve} (\pt \mu^\ve-\overline{\pt \mu^\ve}),\\
		\label{sistem difference 2}
		&M=-\Delta \Phi +\left(\Psi'(\varphi^\ve)-\Psi'(\vphi^\ast ) \right).
	\end{alignat}
\end{subequations}
Testing \eqref{sistem difference 1} by $\Phi-\overline{\Phi}$, we have
\begin{equation}
\label{servo}
(M,\Phi- \overline{\Phi})=-(\G_{\varphi^\ve} (\pt\vphi^\ve-\overline{\pt\vphi^\ve})
-\G_{\vphi^\ast } (\pt \vphi^\ast),\Phi-\overline{\Phi})
-\ve(\G_{\varphi^\ve} (\pt \mu^\ve-\overline{\pt \mu^\ve}),\Phi-\overline\Phi).
\end{equation}
Adding and subtracting $\G_{\vphi^\ve} (\pt \vphi^\ast)$, we can write
$$
(\G_{\varphi^\ve} (\pt \vphi^\ve-\overline{\pt\vphi^\ve})
-\G_{\vphi^\ast} (\pt \vphi^\ast), \Phi- \overline{\Phi})=
(\G_{\vphi^\ve} (\pt \Phi-\overline{\pt\Phi}), \Phi- \overline{\Phi})+((\G_{\varphi^\ve}-\G_{\vphi^\ast}) \pt \vphi^\ast, \Phi- \overline{\Phi}),
$$
where, by the integration by-parts formula \eqref{integration by part formula}, it holds
$$(\G_{\vphi^\ve} (\pt \Phi-\overline{\pt\Phi}),\Phi- \overline{\Phi})=
\frac12\ddt \norm{\sqrt{b(\vphi^\ve)}\nabla \G_{\vphi^\ve} (\Phi- \overline{\Phi})}_{L^2(\Omega)}^2+
\frac12(b'(\vphi^\ve)\pt \vphi^\ve\nabla \G_{\vphi^\ve} (\Phi- \overline{\Phi}),\nabla \G_{\vphi^\ve} (\Phi- \overline{\Phi})).$$
Therefore,  \eqref{servo} is equivalent to
\begin{align}
\label{servo2}
\frac12\ddt \norm{\sqrt{b(\vphi^\ve)}\nabla \G_{\vphi^\ve} (\Phi- \overline{\Phi})}_{L^2(\Omega)}^2+(M,\Phi- \overline{\Phi})
=K_1+K_2+K_3,
\end{align}
having set
\begin{align*}
&K_1=-\frac12(b'(\vphi^\ve)\pt\vphi^\ve\nabla \G_{\vphi^\ve} (\Phi-\overline\Phi),\nabla \G_{\vphi^\ve} (\Phi-\overline\Phi)),\\
 &K_2=-((\G_{\varphi^\ve}-\G_{\vphi^\ast}) \pt \vphi^\ast, \Phi-\overline\Phi),\\
 &K_3=-\ve(\G_{\varphi^\ve} (\pt \mu^\ve-\overline{\pt \mu^\ve}),\Phi-\overline\Phi).
 \end{align*}
Now, multiplying \eqref{sistem difference 2} with $\Phi-\overline\Phi$, by the assumptions on $\Psi$ one obtains
\begin{align*}
	( M,\Phi-\overline\Phi)
	&	=\|\nabla\Phi\|_{L^2(\Omega)}^2
	+( \Psi'(\vphi^\ve)-\Psi'(\vphi^\ast ),\Phi)
	-( \Psi'(\vphi^\ve)-\Psi'(\vphi^\ast),\overline\Phi )
	\\
	&\geq \|\nabla\Phi\|_{L^2(\Omega)}^2 -\theta_0\|\Phi\|_{L^2(\Omega)}^2-|( \Psi'(\vphi^\ve)-\Psi'(\vphi^\ast ),\overline\Phi   ) |\\
	&\geq \|\Phi\|_{H^1(\Omega)}^2 -(\theta_0+1) \|\Phi\|_{L^2(\Omega)}^2-\big(\|\Psi'(\vphi^\ve)\|_{L^1(\Omega)}+\|\Psi'(\vphi^\ast )\|_{L^1(\Omega)}\big)
	|\overline\Phi|.
\end{align*}
Here, we further control $\|\Phi\|_{L^2(\Omega)}^2$ with \eqref{Hilbert interpolation inequality} and \eqref{poincare} as
\begin{align*}
	(\theta_0+1)\|\Phi\|_{L^2(\Omega)}^2\leq C(\|\Phi-\overline{\Phi}\|_{L^2(\Omega)}^2+|\overline\Phi|^2)
	&\leq \frac{1}{4}\|\nabla\Phi\|_{L^2(\Omega)}^2+C\|\Phi-\overline{\Phi}\|_{H^{-1}_{(0)}(\Omega)}^2+C|\overline{\Phi}|^2,
\end{align*}
thus obtaining
\begin{align*}
	( M,\Phi-\overline\Phi)\geq \frac34\|\Phi\|_{H^1(\Omega)}^2- C\|\Phi-\overline{\Phi}\|_{H^{-1}_{(0)}(\Omega)}^2- f(t)|\overline\Phi|,
\end{align*}
where
$$
f(t)=C(1+\|\Psi'(\vphi^\ve(t))\|_{L^1(\Omega)}+\|\Psi'(\vphi^\ast(t))\|_{L^1(\Omega)}).$$
As a result, \eqref{mee} and \eqref{servo2} yield the following differential inequality
\begin{align}
	\label{difference2}
	\frac12\ddt \norm{\sqrt{b(\vphi^\ve)}\nabla \G_{\vphi^\ve} (\Phi-\overline\Phi)}_{L^2(\Omega)}^2+\frac{3}{4}\|\Phi\|_{H^1(\Omega)}^2&\leq
	C\|\Phi-\overline\Phi\|_{H^{-1}_{(0)}(\Omega)}^2+\ve f(t)+K_1+K_2+K_3.
\end{align}
We proceed by estimating the right-hand side in \eqref{difference2}.
To this aim, observe that
$$|K_1|\leq C \int_\Omega |\pt \vphi^\ve||\nabla \G_{\vphi^\ve} (\Phi-\overline\Phi)|^2 \, \d x \le C\|\pt \vphi^\ve\|_{L^6(\Omega)}\|\nabla \G_{\vphi^\ve} (\Phi-\overline\Phi)\|_{L^\frac{12}{5}(\Omega)}^2.$$
An interpolation of $L^\frac{12}{5}(\Omega)$ between $L^2(\Omega)$ and $H^1(\Omega)$ gives
\begin{align*}
\label{R1 estimate}
|K_1|\leq C\|\pt\vphi^\ve\|_{H^1(\Omega)}\|\nabla \G_{\vphi^\ve} (\Phi-\overline\Phi)\|_{L^2(\Omega)}^\frac32
\| \G_{\vphi^\ve} (\Phi-\overline\Phi)\|_{H^2(\Omega)}^\frac12 .
\end{align*}
	Recalling \eqref{equivalence of norm}, we have
	$$\|\G_{\vphi^\ve}(\Phi-\overline\Phi)\|_{H^1(\Omega)}\leq C\|\Phi-\overline\Phi\|_{H^{-1}_{(0)}(\Omega)}.$$
	Besides, by \eqref{H2 norm of G_q-2d} and Theorem \ref{t-uniform} we find
\begin{align*}
\|\G_{\vphi^\ve} (\Phi-\overline\Phi)\|_{H^2(\Omega)}
&\leq C\left( \| \nabla \vphi^\ve\|_{L^2(\Omega)} \|\vphi^\ve\|_{H^2(\Omega)}
\|\nabla\G_{\vphi^\ve} (\Phi-\overline\Phi)\|_{L^2(\Omega)}+\|\Phi-\overline\Phi\|_{L^2(\Omega)}\right) \\
&\leq C \left(\|\nabla\G_{\vphi^\ve} (\Phi-\overline\Phi)\|_{L^2(\Omega)}+\|\Phi-\overline\Phi\|_{L^2(\Omega)}\right)\\
&\leq C\|\Phi-\overline\Phi\|_{H^{-1}_{(0)}(\Omega)}^{\frac12}\|\nabla\Phi\|_{L^2(\Omega)}^{\frac12}.
\end{align*}
Hence, we arrive at
\begin{equation}
	\label{RI}
	K_1\leq C\|\pt \vphi^\ve\|_{H^1(\Omega)} \|\Phi- \overline\Phi\|_{H^{-1}_{(0)}(\Omega)}^{\frac74}\|\nabla\Phi\|_{L^2(\Omega)}^{\frac14}
\leq C\|\pt \vphi^\ve\|_{H^1(\Omega)}^\frac87 \|\Phi- \overline\Phi\|_{H^{-1}_{(0)}(\Omega)}^{2}+\frac18\|\Phi\|_{H^1(\Omega)}^{2}.
\end{equation}
In order to control the term $K_2$, we argue in a similar way as in \eqref{I3 est}. We observe that
$$
K_2=-((\G_{\varphi^\ve}-\G_{\vphi^\ast}) \pt \vphi^\ast,\Phi-\overline\Phi)
=-(\pt \vphi^\ast,(\G_{\varphi^\ve}-\G_{\vphi^\ast }) (\Phi-\overline\Phi))=
((b(\vphi^\ast)-b(\vphi^\ve))\nabla \G_{\varphi^\ve} (\Phi- \overline \Phi),\nabla \mu^\ast).
$$
Indeed, for any $v\in H^1(\Omega)$ with $\overline v=0$, it holds
\begin{align*}
- ( \partial_t \vphi^\ast, &(\G_{\vphi^\ve}-\G_{\vphi^\ast})  v )=
\int_{\Omega} b(\vphi^\ast) \nabla \mu^\ast \cdot \nabla \left( (\G_{\vphi^\ve}-\G_{\vphi^\ast})  v \right) \, \d x
\\
&= \int_\Omega \nabla \mu^\ast \cdot \left( b(\vphi^\ast )-b(\vphi^\ve) \right)  \nabla \G_{\vphi^\ve}   v \, \d x
+
\int_\Omega \nabla \mu^\ast  \cdot b(\vphi^\ve) \nabla \G_{\vphi^\ve} v  \, \d x
-
\int_\Omega \nabla \mu^\ast  \cdot b(\vphi^\ast ) \nabla \G_{\vphi^\ast } v \, \d x
\\
&= - \int_\Omega \nabla \mu^\ast  \cdot \left( b(\vphi^\ve)-b(\vphi^\ast) \right)  \nabla \G_{\vphi^\ve} v \, \d x.
\end{align*}
Then, since $\mu^\ast \in L^\infty(0,T; H^1(\Omega))$ (see \cite[Theorem 2.2]{BB1999}, and \cite[Proposition 4.1]{CGGG2025}), by the
Ladyzhenskaya inequality we have
\begin{align*}
|K_2|&\leq \|\nabla\mu^\ast \|_{L^2(\Omega)} \norm{b(\vphi^\ast )-b(\vphi^\ve)}_{L^4(\Omega)}\norm{\nabla \G_{\varphi^\ve}(\Phi- \overline\Phi)}_{L^4(\Omega)} \\
	&\le C\|\Phi\|_{H^1(\Omega)}\|\G_{\vphi^\ve }(\Phi-\overline\Phi)\|_{H^1(\Omega)}^\frac12
\|\G_{\vphi^\ve }(\Phi-\overline\Phi)\|_{H^2(\Omega)}^\frac12.
\end{align*}
Estimating the norms of $\G_{\vphi^\ve }(\Phi-\overline\Phi)$ as above we obtain
\begin{equation}
	\label{RII}
	K_2\leq C\|\Phi\|_{H^1(\Omega)}^\frac54\|\Phi-\overline\Phi\|_{H^{-1}_{(0)}(\Omega)}^\frac34
\leq C\|\Phi-\overline\Phi\|^2_{H^{-1}_{(0)}(\Omega)}+\frac18 \|\Phi\|^2_{H^1(\Omega)}.
\end{equation}
Finally, we easily control $K_3$ as follows:
\begin{align}
	\label{uno1}
	K_3 \leq \ve \|\nabla\G_{\varphi^\ve}( \pt\mu^\ve-\overline{\pt\mu^\ve})\|_{L^2(\Omega)}\|\Phi-\overline\Phi\|_{H^{-1}_{(0)}(\Omega)}
	\leq C\ve\|\sqrt{\ve}\pt\mu^\ve\|_{L^2(\Omega)}^2+C\|\Phi-\overline\Phi\|_{H^{-1}_{(0)}(\Omega)}^2.
\end{align}
Collecting \eqref{RI}, \eqref{RII}, \eqref{uno1} in the differential inequality \eqref{difference2} and
recalling the equivalence \eqref{equivalence of norm}, we end up with
\begin{equation}
	\label{final diff inq model difference}
	\ddt \norm{\sqrt{b(\vphi^\ve)}\nabla \G_{\varphi^\ve} (\Phi-\overline\Phi)}_{L^2(\Omega)}^2+\|\Phi\|_{H^1(\Omega)}^2\leq
h\norm{\sqrt{b(\vphi^\ve)}\nabla \G_{\varphi^\ve} (\Phi-\overline\Phi)}_{L^2(\Omega)}^2+ \ve f+\ve g,
\end{equation}
where
$$h(t)=C\left(1+\norm{ \pt \vphi^\ve(t)}_{H^1(\Omega)}^\frac87\right)\qquad\text{and}\qquad g(t)=C\norm{\sqrt{\ve}\pt\mu^\ve(t)}_{L^2(\Omega)}^2.$$
Note that $h, f$ and $g$ belong to $L^1(0,T)$, due to Theorem \ref{t-uniform} and \eqref{mutfit2}.
Therefore, an application of the Gr{\"o}nwall Lemma gives
$$\| \Phi(t)-\overline\Phi(t)\| _{H^{-1}_{(0)}(\Omega) }^{2}\leq C \| \Phi (0)-\overline\Phi(0) \| _{H^{-1}_{(0)}(\Omega)}^{2}+
C\ve, \quad \forall \, t\in [0,T].$$
Since $\Phi(0)=0$,
invoking once again \eqref{mee}, we conclude that
$$\| \Phi (t)\| _{(H^1(\Omega))'}\leq C{\sqrt{\ve}}.$$
A subsequent integration of \eqref{final diff inq model difference} on $[0,T]$ gives us the desired inequality and concludes the proof.
$\Box$

\section{Degenerate mobility}
\label{deg mob sec}

In this section, we focus on the degenerate mobility case. Our goal is to establish an existence result for weak solution to system \eqref{w problem}--\eqref{w problem 2}. Some achievements were already obtained in \cite{SLYM}, where the authors considered the system with degenerate mobility and polynomial potential $\Psi$. However, they did not fully exploit the additional regularity guaranteed by the introduction of the $\ve$-perturbation. Herein we will show that it is possible to identify the chemical potential $\mu$, in contrast with what happens in the case of the classical Cahn--Hilliard equation.
\medskip

\textbf{Proof of Theorem \ref{th-degenerate}.} Let $\ve>0$ be fixed. We reason as in \cite{EG1996}, but for the approximating sequence of positive mobilities
$$
b_\lambda(s)=
\begin{cases}
b(-1+\lambda)&s\in[-1,-1+\lambda),\\
b(s)&s\in[-1+\lambda,1-\lambda],\\
b(1-\lambda)&s\in(1-\lambda,1],
\end{cases}
\quad 0<\lambda<1,
$$
extended to $\R$ by $b_\lambda(s)=b(-1+\lambda)$ for $s<-1$ and $b_\lambda(s)=b(1-\lambda)$ for $s>1$, and the regularized potential $\Psi_\lambda(s)=F_\lambda(s)-\frac{\theta_0}{2}s^2$, where $F_\lambda\in C^2(\R)$ is defined by
$$
F_\lambda''(s)=\frac{\theta}{b_\lambda(s)}\quad \forall s\in\R,\qquad F_\lambda(0)=F_\lambda'(0)=0,
$$
so that $F_\lambda=F$ on $[-1+\lambda,1-\lambda]$ and $b_\lambda F''_\lambda\equiv\theta$. 
\smallskip

\textbf{Part 1: $\lambda$-uniform estimates.} 
By the analogue of Theorem \ref{ws thm phi and mu}- (A) for regular potentials, let us consider the corresponding weak solutions $(\vphi_\lambda,\mu_\lambda)$, which satisfy
\begin{align}
	\label{CH1-p lambda}
	&(\pt \vphi_\lambda, \eta) +\l \ve \pt \mu_\lambda , \eta\r + (b_\lambda(\vphi_\lambda)\nabla \mu_\lambda , \nabla \eta)=0 \quad \, \forall \, \eta \in H^1(\Omega), \ \text{a.e. in } (0,\infty),
	\\
	\label{CH2-p lambda}
	&\mu_\lambda = - \Delta \vphi_\lambda + \Psi'_\lambda(\vphi_\lambda)
	\quad \text{a.e. in } \Omega \times (0,\infty).
\end{align}
Reasoning as in the proof of Theorem \ref{ws thm phi and mu} (see, for example Remark \ref{lambda independence 2}), we infer that  the following $\lambda$-independent bounds hold:
\begin{align*}
&\vphi_\lambda \in L^{\infty}(0, \infty; H^2(\Omega)),  \quad \pt \vphi_\lambda \in L^2(0, \infty; H^1(\Omega))\\
&\sqrt{\ve}\mu_\lambda \in L^{\infty}(0, \infty; L^2(\Omega)),\quad \sqrt{b_\lambda(\vphi_\lambda)}\nabla\mu_\lambda \in L^2(0, \infty; L^2(\Omega)),  \quad \pt \mu_\lambda \in L^2(0, \infty; (H^1(\Omega))'), \\
&F_\lambda'(\vphi_\lambda) \in L^{\infty}(0, \infty; L^2(\Omega)).
\end{align*}
Moreover, for any $\lambda>0$, we have the energy balance (cfr. \eqref{"Energy Balance"}):
\begin{equation}
	\label{Energy dissipation degenerate approximation}
E_\lambda(\varphi_\lambda(t), \mu_\lambda(t)) + \int_0^t \left\| \sqrt{b_\lambda(\varphi_\lambda)} \nabla \mu_\lambda \right\|_{L^2(\Omega)}^2 \, \d s
= E_\lambda(\varphi_\lambda(0), \mu_\lambda(0)),
\end{equation}
for any $t\geq 0$, 
where
$$
E_\lambda(\varphi_\lambda, \mu_\lambda) = \frac12\| \nabla \varphi_\lambda  \|_{L^2(\Omega)}^2 + \frac{\varepsilon}{2} \| \mu_\lambda \|_{L^2(\Omega)}^2 + \int_\Omega \Psi_\lambda(\varphi_\lambda) \,\d x.
$$
Indeed, firstly we can reason as in Section \ref{sect-existence}, where we only proved that given any positive time $T>0$ we can find bounds depending on $T$ itself, and then observe that in order to obtain global in time estimates we can follow the approach of Theorem \ref{ws thm phi and mu}.
\begin{remark}
The uniform bound on $F'_\lambda(\vphi_\lambda)$ is a key estimate in order to identify the chemical potential $\mu$ in the next step. This arises from the  $\ve$-perturbation, which enables us to prove a $\lambda$-uniform bound on $\norm{\mu_\lambda}_{L^2(\Omega)}$. Notice that this is not available for the classical Cahn--Hilliard equation with degenerate mobility (see \cite{EG1996}).
\end{remark}
\textbf{Part 2: passage to the limit.} Exploiting the estimates above, as $\lambda \to 0^+$ we have the following convergences, up to a subsequence,
\begin{equation}
	\begin{aligned}
		\label{convergeces degenerte case}
		\vphi_\lambda & \overset{\ast}{\rightharpoonup} \vphi \quad && \text{weak-star in } L^\infty(0, T; H^2(\Omega)), \\
		\vphi_\lambda &\to \vphi &&\text{strongly in } C^0([0, T]; H^{2- \alpha}(\Omega)), \ \alpha \in(0,2], \\
		\pt \vphi_\lambda &\rightharpoonup \pt \vphi && \text{weakly in } L^2(0, T; H^1(\Omega)), \\
		\mu_\lambda & \overset{\ast} \rightharpoonup  \mu &&\text{weak- star in } L^\infty(0, T; L^2(\Omega)), \\
		\pt \mu_\lambda &\rightharpoonup \pt \mu \quad&&\text{weakly in } L^2(0, T; (H^1(\Omega))'), \\
		\mathbf{J}_\lambda:=b_\lambda(\vphi_\lambda)\nabla\mu_\lambda &\rightharpoonup \mathbf{J} \quad&&\text{weakly in } L^2(0, T; L^2(\Omega)).
	\end{aligned}
\end{equation}
Next, exploiting $\norm{F'_\lambda(\vphi_\lambda)}_{L^\infty(0, T; L^2(\Omega))} \le C$ and reasoning as in \cite[Section 4]{DebDet}, we deduce that
\begin{equation}
	\label{phi L infty bound deg prob eps}
	|\vphi(x, t)| < 1 \quad \text{a.e.} \ x \in \Omega, \forall \, t>0 .
\end{equation}
Hence, combining \eqref{convergeces degenerte case} and the weak dominated convergence theorem, we find
$$F'_\lambda(\vphi_\lambda) \rightharpoonup F'(\vphi) \quad \text{weakly in } L^2(0, T; L^2(\Omega)).$$
This is the key point needed in order to pass to the limit in \eqref{CH2-p lambda} and identify $\mu$ as
\begin{equation}
\label{mu identification deg case}
\mu = - \Delta \vphi + \Psi'(\vphi) \quad \text{a.e. in } \Omega \times (0, T).
\end{equation}

Next, observing that, since $b_\lambda$ is non-degenerate, the corresponding solution $\mu_\lambda$ is an element of \linebreak$L^2_{\rm uloc}([0, \infty); H^1(\Omega))$, meaning that we can integrate by parts in \eqref{CH1-p lambda} to obtain
\begin{align}
	\label{phi mu lambda formulation with integartio by parts}
	\l\pt (\vphi_\lambda+\ve  \mu_\lambda), v \r - (\mu_\lambda, \div(b_\lambda(\vphi_\lambda) \nabla v))= 0,
\end{align}
for any $v \in H^2(\Omega)$ such that $\partial_\n v =0$ on $\partial \Omega$. Firstly, we will show that
\begin{equation}
	\label{convergence scalar product mu-divergence}
	\int_0^T (\mu_\lambda, \div(b_\lambda(\vphi_\lambda)\bm{\xi})) \, \eta \, \d t \to \int_0^T (\mu, \div(b(\vphi)\bm{\xi})) \, \eta \, \d t \quad \text{as } \lambda \to 0^+,
\end{equation}
for all $\bm{\xi}\in H^1(\Omega; \mathbb{R}^d)$ such that $\bm{\xi} \cdot \n =0$ on $\partial \Omega$, and for all $\eta \in C_c^\infty(0,T)$. 
By the definition of $b_\lambda$ we have that $b_\lambda \to b$ uniformly on $\mathbb{R}$, hence
\begin{equation}
	\label{mobility convergence lambda}
b_\lambda(\vphi_\lambda) \to b(\vphi)\quad \text{strongly in } L^p(\Omega \times (0,T)), \ \forall \, p \in [2, \infty) \text{ and a.e. in } \Omega\times (0, T).
\end{equation}
Moreover, by \eqref{phi L infty bound deg prob eps}, we also know that
\begin{equation}
\label{b' convergence lambda}
b'_\lambda(\vphi_\lambda) \to b'(\vphi)\quad \text{a.e. in } \Omega\times (0, T).
\end{equation}
We notice that
\begin{align*}
	\int_\Omega \mu_\lambda \div(b_\lambda(\vphi_\lambda) \bm{\xi}) \, \d x 
	&= \int_\Omega \mu_\lambda b_\lambda(\vphi_\lambda) \, \div \bm{\xi} \, \d x + \int_\Omega \mu_\lambda b_\lambda'(\vphi_\lambda) \nabla \vphi_\lambda \cdot \bm{\xi} \, \d x. 
\end{align*}
Exploiting \eqref{convergeces degenerte case} and \eqref{mobility convergence lambda}, $\mu_\lambda b_\lambda(\vphi_\lambda) \rightharpoonup \mu b(\vphi)$ weakly in $L^2(\Omega \times (0,T))$, which implies that
\begin{equation*}
	\int_0^T  \int_\Omega \mu_\lambda b_\lambda(\vphi_\lambda) \, \div \bm{\xi} \, \d x \, \eta \, \d t
	\to \int_0^T \int_\Omega \mu b(\vphi) \, \div \bm{\xi}\,  \d x  \, \eta \, \d t \quad \text{as } \lambda \to 0^+.
\end{equation*}
On the other hand, recalling \eqref{b' convergence lambda}, an application of the dominated convergence theorem yields
\begin{equation}
	\label{b' lambda integral convergence}
	b'_\lambda(\vphi_\lambda) \to  b'(\vphi) \quad \text{strongly in }  L^p(\Omega\times(0,T)),\ \forall\,p<\infty.
\end{equation}
Besides, by  \eqref{convergeces degenerte case}, we have
$$
\nabla\vphi_\lambda \to \nabla\vphi  \quad \text{strongly in }  L^q(\Omega \times (0,T)),\ 3\le q<6.
$$
Hence, 
$$
\int_0^T \int_\Omega \mu_\lambda b_\lambda'(\vphi_\lambda) \nabla \vphi_\lambda \cdot \bm{\xi} \, \d x \, \eta \, \d t 
\to  \int_0^T \int_\Omega \mu b'(\vphi) \nabla \vphi \cdot \bm{\xi} \, \d x \, \eta \, \d t \quad \text{as } \lambda \to 0^+,
$$
which proves \eqref{convergence scalar product mu-divergence}. 
Now, choosing $\bm{\xi}= \nabla v$, with $v \in H^2(\Omega)$ such that $\partial_\n v =0$ on $\partial \Omega$ and recalling \eqref{convergeces degenerte case}, we discover that $(\vphi, \mu)$ satisfies
$$ 
\l\pt (\vphi+\ve  \mu), v \r - (\mu, \div(b(\vphi) \nabla v))= 0,
$$
for any $v \in H^2(\Omega)$ such that $\partial_\n v =0$ on $\partial \Omega$, almost everywhere in $(0,T)$.

Furthermore, we claim that 
\begin{align*}
(\mu, \div(b(\vphi) \bm{\xi} ))
&=-(\Delta \vphi b(\vphi) \div \bm{\xi} + \Delta\vphi b'(\vphi) \nabla \vphi\cdot \bm{\xi} + (b\Psi'')(\vphi)\nabla\vphi \cdot \bm{\xi}),
\end{align*}
for any $\bm{\xi} \in H^1(\Omega; \mathbb{R}^d)$ such that $\bm{\xi} \cdot \n=0$ on $\partial \Omega$, almost everywhere in $(0,T)$. Indeed, observing that
\begin{align*}
\int_0^T \int_\Omega \mu_\lambda \div(b_\lambda(\vphi_\lambda) \bm{\xi}) \, \d x \, \eta \, \d t
&= \int_0^T \int_\Omega -\Delta \vphi_\lambda \div(b_\lambda(\vphi_\lambda) \bm{\xi}) \, \d x \, \eta \, \d t
- \int_0^T \int_\Omega b_\lambda(\vphi_\lambda) \Psi''_\lambda(\vphi_\lambda)  \nabla \vphi_\lambda  \cdot \bm{\xi} \, \d x \, \eta \, \d t
\\
&= 
\int_0^T \int_\Omega -\Delta \vphi_\lambda \div(b_\lambda(\vphi_\lambda) \bm{\xi}) \, \d x \, \eta \, \d t
- \int_0^T \int_\Omega (\theta - \theta_0 b_\lambda(\vphi_\lambda)) \nabla \vphi_\lambda  \cdot \bm{\xi} \, \d x \, \eta \, \d t.
\end{align*}
Therefore, by passing to the limit as above, and arguing by density, we reach the desired claim.

Finally, we can also identify
$$
(\mu , \div(b(\vphi)\bm{\xi})) = -(\mathbf{J}, \bm{\xi}).
$$
This follows from \eqref{convergeces degenerte case}, \eqref{convergence scalar product mu-divergence}, and $(\mu_\lambda,\div(b_\lambda(\vphi_\lambda)\bm{\xi}))=-(\mathbf{J}_\lambda,\bm{\xi})$.
Hence, our solution is also a weak solution in the sense of \cite{EG1996}, i.e.
\begin{equation}
	\label{Degenerate CH with eps 1}
	\l\pt (\vphi+\ve  \mu), v \r + (\mathbf{J}, \nabla v) = 0,
\end{equation}
for any $v \in H^2(\Omega)$ such that  $\partial_\n v =0$ on $\partial \Omega$ (and then, by density, for any $v\in H^1(\Omega)$), with $\mathbf{J}$ identified by
\begin{equation}
	\label{Degenerate CH with eps 2}
	\int_\Omega \mathbf{J} \cdot \bm{\xi} \, \d x = \int_\Omega \Delta \vphi \, \div(b(\vphi) \bm{\xi})+ (b\Psi'')(\vphi)\nabla\vphi \cdot\bm{\xi}\,\d x,
\end{equation}
for all $\bm{\xi}\in H^1(\Omega; \mathbb{R}^d)$ such that $ \bm{\xi} \cdot \n=0$ on $\partial\Omega$.
\medskip

\subsection{Convergence in the degenerate case} 
\label{sub-conv-deg}
Lastly we will show that weak solutions $\left(\vphi^\ve, \mathbf{J}^\ve\right)$  to the generalized Cahn--Hilliard equation converge, on any finite time interval $(0, T)$, to weak solution of the classical Cahn--Hilliard equation, even for degenerate mobility. More precisely we will pass to the limit in \eqref{Degenerate CH with eps 1}--\eqref{Degenerate CH with eps 2} written for $\left(\vphi^\ve, \mathbf{J}^\ve\right)$ with associated $\mu^\ve$.
\smallskip

\textbf{Proof of Theorem \ref{th-degenerate-conv}.}
Thanks to the convergences \eqref{convergeces degenerte case} and the weak lower semicontinuity of the norms, 
we can pass to the limit in the energy equation \eqref{Energy dissipation degenerate approximation} obtaining 
\begin{align*}
 \frac12\| \nabla \varphi^\varepsilon (t) \|_{L^2(\Omega)}^2 &+ \frac{\varepsilon}{2} \| \mu^\varepsilon (t)\|_{L^2(\Omega)}^2 + \int_\Omega \Psi(\varphi^\varepsilon(t)) \,\d x
 + \frac{1}{b_M} \int_0^t\int_\Omega |\mathbf{J}^\ve|^2 \, \d x \,\d s 
\\
&\le \frac12\| \nabla \varphi_0 \|_{L^2(\Omega)}^2 + \frac{\varepsilon}{2} \| \mu_0 \|_{L^2(\Omega)}^2 + \int_\Omega \Psi(\varphi_0) \,\d x,
\end{align*}
for almost every $t \in (0,\infty)$. Here, we have used that 
$|\mathbf{J}_\lambda|^2=b_\lambda^2(\vphi_\lambda) |\nabla\mu_\lambda|^2 \le b_M b_\lambda(\vphi_\lambda) |\nabla\mu_\lambda|^2$.
Since $\varphi_0 \in H^2(\Omega)$ with $ -1 < \varphi_0 < 1$ a.e. in $\Omega$ with $\mu_0 \in L^2(\Omega)$, we deduce that
\begin{equation}
	\label{Energy inequality approximation degenerate}
\norm{\nabla \vphi^\ve}_{L^\infty(0,\infty; L^2(\Omega))} 
+\sqrt{\ve} \norm{\mu^\ve}_{L^\infty(0,\infty; L^2(\Omega))} 
+ \norm{\mathbf{J}^\ve}_{L^2(0,\infty;L^2(\Omega))} \le C,
\end{equation}
for some $C>0$ independent of $\varepsilon$. As a consequence, reasoning as in \eqref{media phi}--\eqref{media},  it is easily seen that
$$
|\overline{ \varphi^\varepsilon(t)} - \overline{\varphi_0} | \le C \sqrt{\varepsilon},
\quad \forall \, t \geq 0,
$$
which implies that 
 $\overline{\vphi^\ve}\to \overline{\vphi_0}$
uniformly in $[0,\infty)$ as $\ve\to  0$, and there exist $\widetilde\ve>0$ and $\widetilde{\delta}>0$ such that
\begin{equation}
\label{media-deg}
	\overline{\vphi^\ve(t)} \in [-1+\widetilde{\delta},1-\widetilde{\delta}],\qquad \forall \, t\geq 0, \ \forall \, \ve\in (0,\widetilde\ve).
\end{equation}
Also, by comparison in \eqref{Degenerate CH with eps 1}, 
we infer
\begin{equation}
\label{der-time}
\| \partial_t \varphi^\varepsilon + \varepsilon \partial_t \mu^\varepsilon \|_{L^2(0, \infty; H^{-1}_{(0)}(\Omega))} 
\le 
\| \mathbf{J}^\varepsilon \|_{L^2(0, \infty; L^2(\Omega))} 
\le C,
\end{equation}
for some $C>0$ independent of $\varepsilon$.

In order to pass to the limit in \eqref{Degenerate CH with eps 2} for $\mathbf{J}^\ve$, we need to recover an uniform estimate on $\Delta \vphi^\ve$ in $L^2(0,T; L^2(\Omega))$. As already mentioned in the Introduction, we will propose a new argument based on the dissipation term in the energy equality.
To this aim, let $\left(\vphi_\lambda^\ve, \mu_\lambda^\ve\right)$ be the solution to \eqref{CH1-p lambda}--\eqref{CH2-p lambda}, which converges to 
$\left(\vphi^\ve, \mu^\ve\right)$ as $\lambda \to 0$. In light of \eqref{Energy dissipation degenerate approximation}, we notice that there exists $C>0$ independent of $\ve$ and $\lambda$ such that
\begin{equation}
\label{7.15}
\norm{\nabla \vphi_\lambda^\ve}_{L^\infty(0,\infty; L^2(\Omega))} +
\left\| \sqrt{b_\lambda(\varphi_\lambda^\varepsilon)} \nabla \mu_\lambda^\varepsilon \right\|_{L^2(0, \infty, L^2(\Omega))} \le C.
\end{equation} 
Exploiting $\mu_\lambda^\varepsilon = -\Delta \varphi_\lambda^\varepsilon + \Psi_\lambda'(\varphi_\lambda^\varepsilon)$ and recalling $\Psi''_\lambda(\vphi_\lambda^\ve)= F''_\lambda(\vphi_\lambda^\ve)-\theta_0$ and $F''_\lambda(\vphi_\lambda^\ve)= \frac{\theta}{b_\lambda(\vphi_\lambda^\ve)}$, we infer
\begin{align}
&\int_t^{t+1} \int_\Omega b_\lambda(\vphi_\lambda^\ve)|\nabla \mu_\lambda^\ve|^2 \, \d x \, \d s 
= \int_t^{t+1} \int_\Omega b_\lambda(\vphi_\lambda^\ve)|\nabla ( - \Delta \vphi_\lambda^\ve + \Psi'_{\lambda}(\vphi_\lambda^\ve))|^2 \, \d x \, \d s \nonumber\\
& \quad = \int_t^{t+1} \int_\Omega b_\lambda(\vphi_\lambda^\ve) |\nabla \Delta \vphi_\lambda^\ve|^2 \, \d x \, \d s
 + \int_t^{t+1} \int_\Omega b_\lambda(\vphi_\lambda^\ve) \Psi''_\lambda (\vphi_\lambda^\ve)^2 |\nabla\vphi_\lambda^\ve|^2 \, \d x \, \d s 
 \notag 
 \\
 &\qquad -2\int_t^{t+1} \int_\Omega b_\lambda(\vphi_\lambda^\ve) \Psi''_\lambda (\vphi_\lambda^\ve) \nabla \Delta\vphi_\lambda^\ve \cdot \nabla\vphi_\lambda^\ve \, \d x \, \d s \nonumber\\
	&\quad= \int_t^{t+1}\int_\Omega b_\lambda(\vphi_\lambda^\ve) |\nabla \Delta \vphi_\lambda^\ve|^2 \, \d x \, \d s + \theta \int_t^{t+1} \int_\Omega F''_\lambda(\vphi_\lambda^\ve) |\nabla \vphi_\lambda^\ve|^2 \, \d x \, \d s
	- 2 \theta \theta_0 \int_t^{t+1} \int_\Omega |\nabla \vphi_\lambda^\ve|^2 \, \d x \, \d s 
	 \nonumber \\
	& \qquad + \theta_0^2 \int_t^{t+1} \int_\Omega b_\lambda(\vphi_\lambda^\ve)|\nabla\vphi_\lambda^\ve|^2 \, \d x \, \d s 
	- 2 \theta \int_t^{t+1} \int_\Omega \nabla \Delta \vphi_\lambda^\ve \cdot \nabla \vphi_\lambda^\ve \, \d x \, \d s \nonumber \\
	&\qquad + 2 \theta_0 \int_t^{t+1} \int_\Omega b_\lambda(\vphi_\lambda^\ve) \nabla\Delta \vphi_\lambda^\ve \cdot  \nabla \vphi_\lambda^\ve \, \d x \, \d s.
	\label{conto smart}
\end{align}
We highlight that the above computations are rigorous since $\vphi^\ve_\lambda \in L^2(0,T; H^3(\Omega))$, which follows from the elliptic regularity applied to  \eqref{CH2-p lambda}.
Then, an integration by parts yields
\begin{equation}
\label{stima-deg}
- 2 \theta \int_t^{t+1} \int_\Omega \nabla \Delta \vphi_\lambda^\ve \cdot \nabla \vphi_\lambda^\ve \, \d x \, \d s
=
2 \theta\int_t^{t+1} \int_\Omega |\Delta \vphi_\lambda^\ve|^2 \, \d x \, \d s. 
\end{equation}
On the other hand, by the Young inequality, we have
\begin{equation*}
2 \theta_0  \left| \int_t^{t+1} \int_\Omega b_\lambda(\vphi_\lambda^\ve) \nabla\Delta \vphi_\lambda^\ve \cdot  \nabla \vphi_\lambda^\ve \, \d x \, \d s\right|
\leq 
\frac12 \int_t^{t+1} \int_\Omega b_\lambda(\vphi_\lambda^\ve) |\nabla \Delta \vphi_\lambda^\ve|^2 \d x \, \d s + 2 \theta_0^2 b_M \int_t^{t+1} \int_\Omega |\nabla \vphi_\lambda^\ve|^2\,  \d x \, \d s.
\end{equation*}
Thus, we arrive at
\begin{align*}
&
2 \theta \int_t^{t+1} \int_\Omega |\Delta \vphi_\lambda^\ve|^2 \, \d x \, \d s +
\frac12 \int_t^{t+1}\int_\Omega b_\lambda(\vphi_\lambda^\ve) |\nabla \Delta \vphi_\lambda^\ve|^2\, \d x \, \d s 
+  \int_t^{t+1} \int_\Omega \left( \theta F''_\lambda(\vphi_\lambda^\ve) +\theta_0^2 b_\lambda(\vphi_\lambda^\ve) \right)  |\nabla \vphi_\lambda^\ve|^2 \, \d x \, \d s \nonumber \\
& \leq 	\int_t^{t+1} \int_\Omega b_\lambda(\vphi_\lambda^\ve) |\nabla \mu_\lambda^\ve|^2 \,  \d x \, \d s + 2 \theta_0(\theta+ \theta_0b_M) \int_t^{t+1} \int_\Omega|\nabla \vphi_\lambda^\ve|^2\,  \d x \, \d s. 
\end{align*}
Since $F''_\lambda (s) \ge 0$, and by using \eqref{7.15},
we obtain
\begin{equation*}
	\norm{\Delta\vphi_\lambda^\ve}_{L^2_{\rm uloc}([0, \infty); L^2(\Omega))} \le C,
\end{equation*}
where $C$ is independent of $\lambda$ and $\varepsilon$. Passing to the limit in $\lambda$, and exploiting \eqref{Energy inequality approximation degenerate} and \eqref{media-deg}, we infer
\begin{equation}
	\label{unifor estimates phi degenerate problem}
	\norm{\vphi^\ve}_{L^\infty(0, \infty; H^1(\Omega))}+ \norm{\vphi^\ve}_{L^2_{\rm uloc}([0, \infty); H^2(\Omega))} \le C.
\end{equation}

Now, by using \eqref{Energy inequality approximation degenerate}, \eqref{der-time} and \eqref{unifor estimates phi degenerate problem}, we deduce, up to subsequence, the following convergences for any $T>0$:
\begin{equation}
\label{weak convergences degenerate case}
\begin{aligned}
\varepsilon \mu^\varepsilon &\to 0 \quad &&\text{strongly in } L^\infty(0, T; L^2(\Omega)) \cap C([0,T];(H^1(\Omega))'),
\\
\vphi^\ve &\overset{\ast}{\rightharpoonup} \vphi^\ast  && \text{weak-star in } L^\infty(0, T; H^1(\Omega)),\\
	\vphi^\ve &\rightharpoonup \vphi^\ast  && \text{weakly in } L^2(0, T; H^2(\Omega)),\\
	\mathbf{J}^\ve &\rightharpoonup \mathbf{J}^\ast  && \text{weakly in } L^2(0, T; L^2(\Omega)),\\
	\pt(\vphi^\ve+ \ve \mu^\ve) &\rightharpoonup \pt \vphi^\ast  \quad&& \text{weakly in } L^2(0, T; (H^1(\Omega))').
	\end{aligned}
\end{equation}
For the latter convergence, we reasoned as in \eqref{phi limit has a time derivative}. Moreover, by the Aubin-Lions lemma, $\vphi^\ve+ \ve \mu^\ve \to\vphi^\ast$ in $C([0,T];(H^1(\Omega))')$. 
Then, 
we also deduce that $\vphi^\ve\to\vphi^\ast$ in $C([0,T];(H^1(\Omega))')$. Arguing as in Section \ref{S-convergenza}, 
we infer
\begin{equation}
	\label{phi strong convergence deg case}
	\vphi^\ve \to \vphi^\ast  \quad \text{strongly in }  C([0, T]; L^2(\Omega))\cap L^2(0, T; H^{2- \alpha}(\Omega)), \ \forall \, \alpha \in (0,2].
\end{equation}
Thus $\vphi^\ve \to \vphi^\ast$ a.e. in $\Omega \times (0,T)$ which implies $\vphi^\ast \in [-1,1]$. Since $b \in C^1([-1,1])$, this entails that
\begin{equation}
	\label{convergence of mobility degenerate case}
	b(\vphi^\ve) \to b(\vphi^\ast) \quad \text{and} \quad b'(\vphi^\ve) \to b'(\vphi^\ast) \quad \text{a.e. in } \Omega \times (0, T).
\end{equation}
As a consequence, we can pass to the limit in \eqref{Degenerate CH with eps 1} to discover that $\left(\vphi^\ast, \mathbf{J}^\ast \right)$ solves
\begin{equation}
	\int_0^T \l \pt \vphi^\ast , v\r \, \d s + \int_0^T \int_\Omega \mathbf{J}^\ast  \cdot \nabla v \,  \d x \, \d s =0, \quad \forall \, v \in L^2(0, T; H^1(\Omega)).
\end{equation}
Next, we identify $\mathbf{J}^\ast $ passing to the limit in (cf. \eqref{Degenerate CH with eps 2})
\begin{align}
	\label{J eps}
	\int_0^T \int_\Omega \mathbf{J}^\ve \cdot \bm{\xi} \, \d x \, \d s &=  \int_0^T \int_\Omega \Delta \vphi^\ve \, b(\vphi^\ve)  \div \bm{\xi} + \Delta \vphi^\ve \, b'(\vphi^\ve) \nabla \vphi^\ve \cdot \bm{\xi}  + (b\Psi'')(\vphi^\ve)\, \nabla \vphi^\ve \cdot \bm{\xi} \, \d x \, \d s\nonumber\\
	&= Z_1+ Z_2 + Z_3,
\end{align}
for any $\bm{\xi} \in L^2(0, T; H^1(\Omega;\mathbb{R}^d)) \cap L^\infty(\Omega \times (0, T); \mathbb{R}^d)$ such that $ \bm{\xi} \cdot \n  =0 $ on $\partial \Omega \times (0,T)$.

Firstly, concerning $Z_1$, by the dominated convergence theorem 
we know that $b(\vphi^\ve) \div \bm{\xi} \to b(\vphi^\ast) \div \bm{\xi}$ strongly in $L^2(0,T;L^2(\Omega))$. Then, by \eqref{weak convergences degenerate case}, it follows that
\begin{equation}
Z_1 \to \int_0^T \int_\Omega \Delta \vphi^\ast \,  b(\vphi^\ast ) \, \div \bm{\xi} \, \d x \, \d s.
\end{equation}
Regarding $Z_2$, since $\nabla \vphi^\ve \to \nabla \vphi^\ast$ strongly in $L^2(0,T; L^2(\Omega))$ by \eqref{phi strong convergence deg case}, and recalling \eqref{convergence of mobility degenerate case}, the dominated convergence theorem entails that
\begin{align}
\int_0^T \int_\Omega |b'(\vphi^\ve)\nabla \vphi^\ve 
- b'(\vphi^\ast )\nabla \vphi^\ast |^2 \, \d x \, \d s
\to 0, \quad \text{as } \ve \to 0.
\end{align}
Hence, since $\bm{\xi} \in L^\infty(\Omega \times (0, T); \mathbb{R}^d)$,
by using \eqref{weak convergences degenerate case}, we get
\begin{equation*}
Z_2= \int_0^T\int_\Omega \Delta\vphi^\ve \, b'(\vphi^\ve) \nabla\vphi^\ve \cdot \bm{\xi} \, \d x \, \d s \to 
\int_0^T\int_\Omega \Delta\vphi^\ast \, b'(\vphi^\ast ) \nabla\vphi^\ast \cdot  \bm{\xi} \, \d x \, \d s\quad \text{as } \ve \to 0.
\end{equation*}
Concerning $Z_3$, we observe that, since $|\vphi^\ve(x, t)|<1$ a.e. in $\Omega\times(0, T)$,
$$
\int_0^T\int_\Omega (b\Psi'')(\vphi^\ve)  \nabla \vphi^\ve \cdot \bm{\xi} \, \d x \, \d s
 = \int_0^T \int_\Omega (\theta-\theta_0b(\vphi^\ve))\nabla\vphi^\ve \cdot \bm{\xi} \, \d x \, \d s.
 $$ 
Thus, \eqref{phi strong convergence deg case} and \eqref{convergence of mobility degenerate case} are sufficient to infer
\begin{equation*}
Z_3= \int_0^T\int_\Omega (b \Psi'')(\vphi^\ve) \nabla \vphi^\ve \cdot\bm{\xi} 
	\, \d x \, \d s \to 
	\int_0^T \int_\Omega(\theta-\theta_0b(\vphi^\ast)) 
	\nabla\vphi^\ast \cdot\bm{\xi} \, \d x \, \d s
	=\int_0^T\int_\Omega (b \Psi'')(\vphi^\ast) \nabla \vphi^\ast \cdot \bm{\xi} \, \d x \, \d s.
\end{equation*}
In conclusion it holds
\begin{align*}
 \int_0^T \int_\Omega \mathbf{J}^\ast \cdot \bm{\xi} \, \d x \, \d s 
&=
\int_0^T \int_\Omega \Delta \vphi^\ast \, b(\vphi^\ast) \, \div \bm{\xi} \,  \d x \, \d s + \int_0^T\int_\Omega \Delta\vphi^\ast \,
b'(\vphi^\ast) \nabla \vphi^\ast \cdot \bm{\xi} \, \d x \, \d s 
\\
&
\quad + \int_0^T\int_\Omega (b \Psi'')(\vphi^\ast) \nabla \vphi^\ast \cdot \bm{\xi} \, \d x \, \d s,
\end{align*}
for any $\bm{\xi} \in L^2(0, T; H^1(\Omega;\mathbb{R}^d)) \cap L^\infty(\Omega \times (0, T);\mathbb{R}^d)$ such that $ \bm{\xi} \cdot \n  =0 $ on $\partial \Omega\times(0, T)$. Therefore, $\left(\vphi^\ast, \mathbf{J}^\ast \right)$ is a weak solution of the classical Cahn--Hilliard equation with degenerate mobility in the sense of \cite{EG1996}.

\bigskip

\noindent
\textbf{Acknowledgments.}
This work is partially supported by the GNAMPA (Gruppo Nazionale per l'Analisi Ma\-te\-ma\-ti\-ca, la Probabilit\`{a} e le loro Applicazioni) of INdAM (Istituto Nazionale di Alta Matematica) project CUP E5324001950001. A. Giorgini and M. Conti are supported by the MUR grant Dipartimento di Eccellenza 2023-2027 of Dipartimento di Matematica, Politecnico di Milano.

\end{document}